\documentclass[a4paper,12pt]{amsart}

\usepackage[pctex32]{graphics}
\usepackage{graphics}
\usepackage[dvips]{geometry}
\usepackage{psfrag}
\usepackage{graphicx}
\usepackage{amscd}
\usepackage{amsmath,amstext, amsthm}
\usepackage{amsfonts}
\usepackage{amssymb}
\usepackage{epsfig} 

\newtheorem{maintheorem}{Theorem}
 
\newtheorem{theorem}{Theorem}[section]
\newtheorem{corollary}[theorem]{Corollary}

\newtheorem{lemma}[theorem]{Lemma}
\newtheorem{proposition}[theorem]{Proposition}

\theoremstyle{definition}
\newtheorem{definition}[theorem]{Definition}
\newtheorem{remark}[theorem]{Remark}

\numberwithin{equation}{section}

\newcommand{\vfi}{\varphi}

\newcommand{\ga}{\gamma}
\newcommand{\la}{\lambda}

\newcommand{\de}{\delta}

\renewcommand{\epsilon}{\varepsilon}

\newcommand{\const}{\operatorname{const}}

\newcommand{\F}{\mathcal{F}}
\newcommand{\U}{\mathcal{U}}

\newcommand{\Le}{\mathcal{L}}
\newcommand{\R}{I\!\!R}

\newcommand{\D}{\mathcal{D}}
\newcommand{\C}{\mathcal{C}}

\newcommand{\K}{\mathcal{K}}
\newcommand{\B}{\mathcal{B}}

\newcommand{\V}{\mathcal{V}}
\newcommand{\T}{\mathcal{T}}

\begin{document}

\title{Mixing rate for semi-dispersing billiards with non-compact
 cusps}
\author{A. Arbieto, R. Markarian, M. J. Pacifico \and R. Soares}


\thanks{This work was partially supported by CNPq Brazil, Pronex on Dynamical Systems, 
FAPERJ-Cientista do Nosso Estado, E-26/100.588/2007, FAPERJ-Bolsa Nota 10-Proc E-26/100.053/2007,
Prodoc-CAPES and Proyecto PDT 2006-2008. S/C/IF/54/001, Uruguay.}


\begin{abstract}   
Since the seminal work of Sinai one studies chaotic
properties of planar billiards tables. Among them  is the study of
 decay of correlations for these tables. There are examples 
in the literature of tables 
with exponential and even polynomial decay.

However, until now nothing is known about mixing properties
for billiard tables with non-compact cusps. There is no consensual
definition of mixing for systems with infinite 
invariant measure. In this paper we study geometric and ergodic 
properties of billiard tables with a 
non-compact cusp. 
The goal of this text is, using the definition of mixing proposed
by Krengel and Sucheston for systems with invariant infinite measure,
to show that the billiard whose table is constituted by the
$x$-axis and and the portion in the plane below the graph
of $f(x)=\frac{1}{x+1}$ is mixing and the speed of mixing is
polynomial.
\end{abstract}

\maketitle


\medskip

\section{Introduction} \label{sec:introduction} \pagenumbering{arabic}

The \emph{planar billiard} is the dynamical system defined
by the free motion of a particle in the interior of a domain
$\mathcal{D} \subset
\R^2$ (usually called {\em{table}}) subjected to elastic collisions to the boundary of $\D$,
 that is, angle of incidence equals angle of reflexion.
In a seminal work, Sinai \cite{Si70} proved that the billiard map
of a system in a two-dimensional torus with finitely many convex
obstacles is a K-automorphism. 

For billiards with non-compact cusps, that generate a dynamical system
with an infinite invariant measure, in \cite{Le02} Lenci proved an
extension of the results of Katok and Strelcyn \cite{Ka86}
for the infinite measure case and, as an application, he showed
that certain tables with non-compact cusps have hyperbolic
structure, that is, existence of absolutely
continuous local stable and unstable manifolds. Furthermore,
adapting arguments contained in \cite{LW95}, Lenci proved that
these billiards maps are ergodic.

About the finite measure case, in \cite{BS81},
Bunimovich and Sinai proved a ``stretched'' exponential decay of correlations
for dispersing billiards.
Young \cite{Yo98} showed that the decay of correlations is actually exponential.
This later result was extended by Chernov \cite{Ch99} for billiards with positive-angle
corners. 

In \cite{Ma04}, Markarian, based on \cite{Yo99}, showed that 
 billiards in the Bunimovich stadium has polynomial decay of correlations. 
 More recently, Chernov and Markarian \cite{CM07} proved that semi-dispersing 
billiard tables with compact cusps also have polynomial decay 
 of correlations. 
Improved estimates for correlations in different types of billiard
tables 
were also proved by  Chernov and Zhang \cite{CZ08}.

%
%
%
%

We are interested in tables of the form $\D = \{(x,y) \in
\R^2 : x \geq 0, 0 \leq y \leq f(x)\},$ where $f: \R_0^+
\rightarrow \R^+$ is a three times differentiable bounded convex function,
satisfying the hypotheses (H1) to (H5) listed in Section \ref{sec:definition}.

\begin{maintheorem}  \label{main:ksystem}
 The billiard map defined in a table $\D$ with a non-compact cusp is
an infinite K-automorphism.
\end{maintheorem}

Following Krengel and Sucheston \cite{KS69}, we say that an endomorphism $\F$ 
on a $\sigma$-finite infinite measure space $(X,\B,\mu)$ is \emph{F-mixing}
if for all measurable set $A$ with $\mu(A) < \infty$,
$$\mu(\F^n A \cap A) \to 0 \quad \mbox{as} \quad n \to \infty.$$ 

As it was commented before, there is no consensual  definition of
mixing for systems with infinite measure. A discussion on different 
definitions of mixing for systems with infinite measure was recently 
done by Lenci \cite{Le09}.

Parry \cite{Pa65} showed that \emph{``an infinite K-automorphism has
countable Lebesgue spectrum"} and in \cite{KS69}, Krengel and Sucheston
showed that \emph{``if an endomorphism has countable Lebesgue spectrum then this endomorphism is
F-mixing''}. Therefore we get 

\begin{corollary} \label{main:mixing}
The billiard map defined in a table $\D$ with a non-compact cusp is F-mixing.
\end{corollary}

 For a conservative endomorphism $\F$ on a $\sigma$-finite measure space $(X,\B,\mu)$, we define its 
{\it entropy} \cite{Kr67} by
\[
 h(\F) = \sup\{h(\F_E,\mu_E) \;|\;\; E \subset X, 0 < \mu(E) <\infty\}.
\]
In \cite[p. 172]{Kr67}, Krengel showed that \emph{``every conservative K-automorphism
on a $\sigma$-finite measure space has positive entropy"}. Hence

\begin{corollary}
The entropy of the billiard map defined in a table $\D$ with a non-compact cusp
is positive.
\end{corollary}

Furthermore we can study the speed of convergence to zero in this definition of
F-mixing. We say that an endomorphism $\F$ is {\it polynomially F-mixing} if
 \[
\mu(\F^{-n}A \cap A) \geq C\frac{1}{n^\alpha},
\]
for some ``good" (e.g., with piecewise differentiable boundary)
set $A$ with $0<\mu(A) < \infty$ and some  
$\alpha >0$. The constant $C$  
depends on $A$ but the exponent $\alpha$ depends only on $\F$. 

Using $f(x)= (x+1)^{-1}$ in the definition of the table $\D$ we show

\begin{maintheorem} \label{main:decay}
The billiard map defined in a table $\D$ with a non-compact cusp is polynomially F-mixing.
\end{maintheorem}

\section{Definition of the dynamical system} \label{sec:definition}

As mentioned in the previous section, we are interested in tables of the form $\D = \{(x,y) \in
\R^2 : x \geq 0, 0 \leq y \leq f(x)\},$ where $f: \R_0^+
\rightarrow \R^+$ is a three times differentiable bounded convex function.

\begin{figure}[ht]
\begin{center}
\includegraphics[height=2.0in]{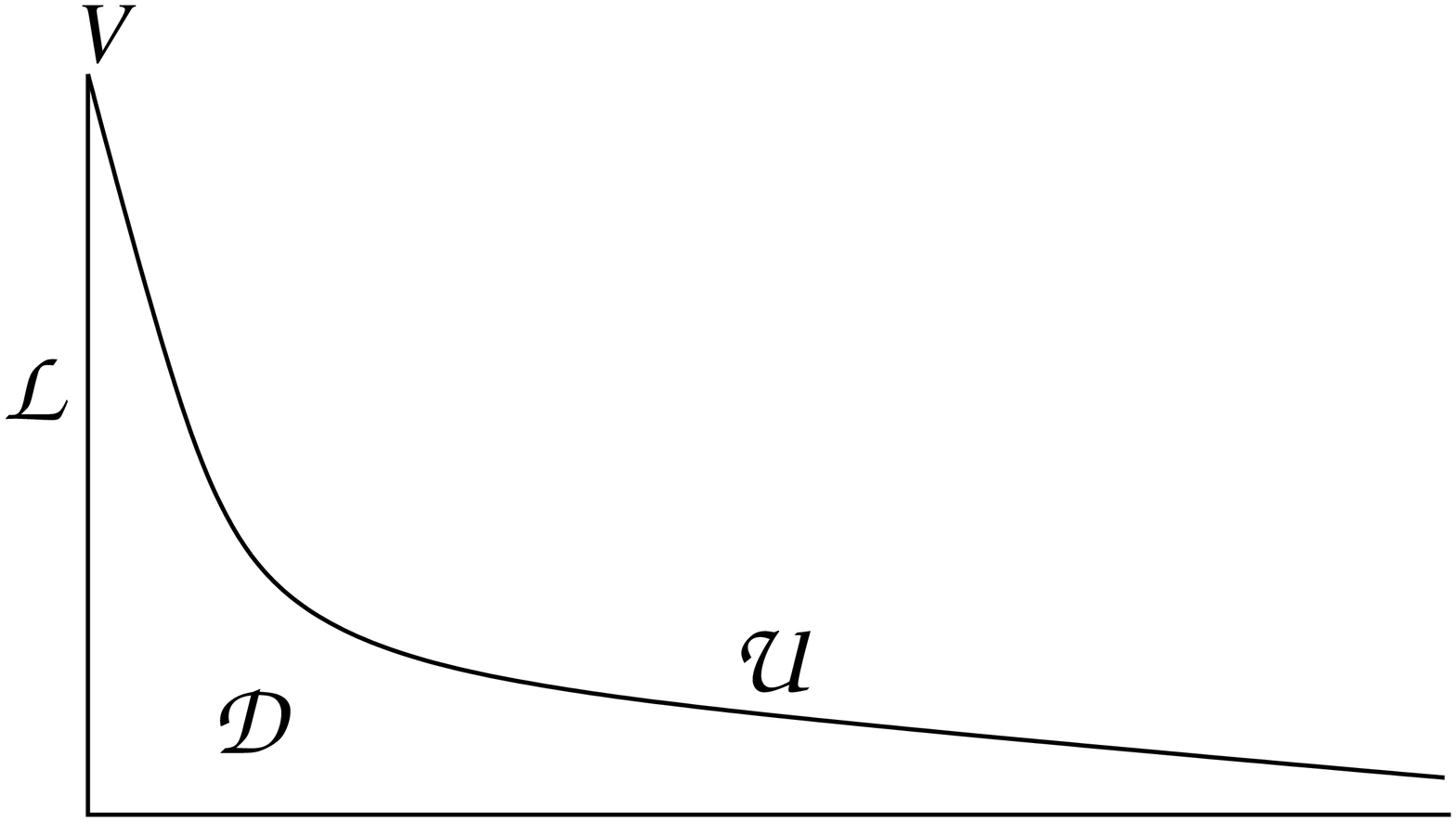}
\caption{\label{f:q} Introducing $\D$, $\U$, $\Le$ e $V$.}
\end{center}
\end{figure}

We denote by $\U$ the dispersing part of the table $\D$ and by $\Le$ the leftmost vertical
wall in $\D$.
The angle in the vertex $V=(0,f(0))$ is $\pi/2+\arctan f'(0^{+})$ and it can be zero.
So the billiard table might have a compact cusp besides the non-compact one on
$x=+\infty$.

 We present two other tables, that will be used in the definitions below: 
                 \[\D_{2}=\{(x,y)\in\R^{+}_{0} \times \R:|y|\leq f(x)\},\]
                 \[\D_{4}=\{(x,y)\in\R \times \R :|y|\leq f(|x|)\}.\]

%

For $f, g: \R_0^+ \rightarrow \R^+ $ we use the following notations: $f(x)<\!<g(x)$ indicates
that there exists a constant $C$ such that $f(x) \leq Cg(x)$, as $ x\to \infty$,
analogously for the symbol $>\!>$ and we denote by  $f = o(g)$ if $\frac{f(x)}{g(x)}$
tends to zero, as $x \to \infty$. Moreover, we use the same symbols when $x \to 0$,
if there is no ambiguity.
Also, we indicate by $A \asymp B$ if there exists a constant
 $C>0$ such that $C^{-1} < A/B < C$ and
we write $A= O(B)$ if there exists a constant $C>0$ such that $|A|/B <C$.


Define $x_{t}=x_{t}(x)$, for each $x$ on $\D_{2}$, implicitly by
\[
\frac{f(x)+f(x_{t})}{x-x_{t}}=-f'(x_t).
\]
One can see that $x_t$ is the $x$-coordinate of the tangent point on $\U$.	

\begin{figure}[ht]
\begin{center}
\includegraphics[scale=0.5]{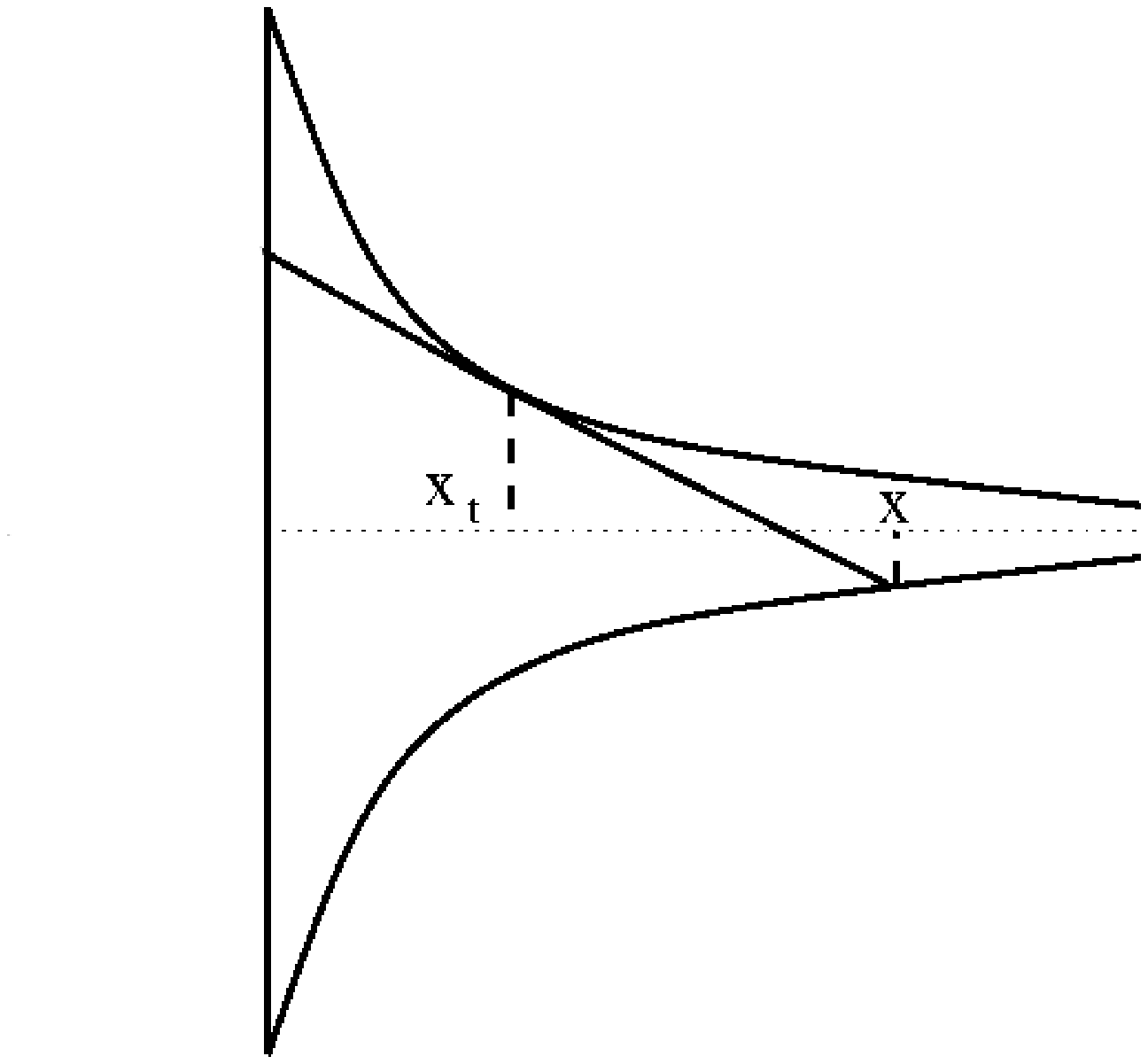}
\caption{\label{fig:xu} The point $x_t$.}
\end{center}
\end{figure}

In \cite{Le02}, Lenci studied tables with $f: \R_0^+ \rightarrow \R^+$
satisfying the following assumptions
\begin{enumerate}
\item[(H1)] $f^{\prime\prime}(x) \to 0$ as $x \to +\infty$;
 \item[(H2)] $|f^{\prime}(x_t)| <\!<
|f^{\prime}(x)|$;
 \item[(H3)]
$\dfrac{f(x)f^{\prime\prime}(x)}{(f^{\prime}(x))^2} >\!> 1$;
\item[(H4)]
$\dfrac{|f^{\prime\prime\prime}(x)|}{f^{\prime\prime}(x)} <\!< 1$;
\item[(H5)] $|f^{\prime}(x)| >\!> (f(x))^{\theta}$, for some
$\theta
> 0$.
\end{enumerate}
It is not difficult to see that $\displaystyle f(x)=\frac{1}{x+1}$ satisfies the conditions
above.

\begin{figure}[ht]
\begin{center}
\includegraphics[height=1.0in]{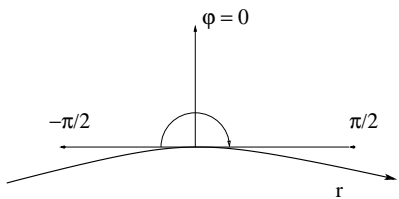}
\caption{\label{fig:orientation} the choice of orientation for $r$ and $\varphi$}
\end{center}
\end{figure}

Following \cite{Le02}, choosing as cross-section the rebounds against the dispersing part
$\U$. we parametrize these line elements as 
 $z = (r, \varphi)$, $r \in (-\infty,0]$ is the arc length variable along $\U$ (with $r=0$ for 
 the vertex $V$) 
and $\varphi \in [-\pi/2,\pi/2]$ is the angle between the velocity vector and the normal
at the point of collision, as in Figure \ref{fig:orientation}. We define the manifold
$M = (-\infty,0)\times (-\pi/2,\pi/2)$ and the return map $T$ defined on $M$,
preserving the measure $d\mu = \cos \varphi drd\varphi$.

\begin{figure}[ht]
\begin{center}
\includegraphics[width=4.0in,height=2.0in]{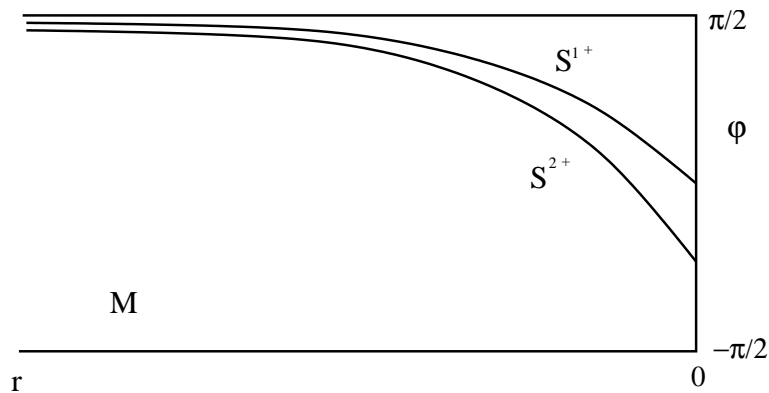}
\caption{\label{f:sing} Singularity lines.}
\end{center}
\end{figure}

We do not define $T$ on those points that hit tangentially $\U$ or that would end
up in the vertex $V$. That is, we exclude $T^{-1}\partial M$. These points make up
the singularity set of $T$, denoted by $S$. This set consists of two lines
(see \cite[p.138]{Le02}) $S^{+}=S^{1+}\cup\; S^{2+}$ (as shown in Figure \ref{f:sing}).
 The curve $S^{1+}$ corresponds to tangencies on $\partial \D_4$ in the third
quadrant (on $\D$, tangencies on $\U$, after a rebound on the vertical side); 
this curve is as regular as $f$. As for $S^{2+}$, its first part corresponds
to line elements pointing to $V$ (on $\D$, after a rebound on the horizontal side); 
as $r$ decreases, these become tangencies on $\partial \D$. the boundary between these two behaviors
is the only non-regular point of $S^{2+}$. 

Analogously we define $S^{-}=S^{1-}\cup\; S^{2-}$, where $S^{i-}, i=1,2$ are the singularity lines
of $T^{-1}$, obtained from $S^{i+}$ using the time-reversal operator
$(r, \varphi) \mapsto (r, -\varphi)$.
We denote by $S^{\pm}_{n} = \bigcup_{i=0}^n T^{\mp
i}S^{\pm}$ e $S^{\pm}_{\infty} = \lim_{n \to \infty} S^{\pm}_{n}$.

%

On $TM$ we define the cone bundles \cite[Section 4]{Le02},
\begin{center}
$\C^u(z) =\{(dr,d\varphi) \in T_zM:drd\varphi \geq 0\}$\\
$\C^s(z) = \{(dr,d\varphi) \in T_zM:drd\varphi \leq 0\}$,
\end{center}
which will be denoted {\it unstable} and {\it stable} cones, respectively.
They are strictly invariant under the action of $T$.

\begin{remark}
We note that our choice of parametrization is different from the one in \cite{Le02}.
This leads to a different choice of the cone bundles.
However, it does not alter the results obtained in that paper.
\end{remark}

Let $\Le$ be the leftmost wall on $\D$ and $M_4$ the phase space defined by
the vectors based on $\Le$. Since $\Le$ is a global cross-section we can define
 a return map $T_{4}$ and let $\mu_{M_4}$ be the measure $\mu$
induced on $M_4$. Denote by $M_3$ the region of $M$ located above $S^{2+}$.
From the definition of $S^{2+}$, the line elements of $M_3$ are precisely
the ones that,on $\D_2$, hits $y$-axis. We call $T_3$
the return map to $M_3$ and one can see that $(M_3, T_3, \mu_{M_3})$
is isomorphic (with respect to $\mu$) to $(M_4, T_{4}, \mu_{M_4})$.

Lenci \cite{Le02} showed that the billiard map $T$ has a 
hyperbolic structure,
 i.e., existence of local stable and unstable
manifolds almost everywhere and these local foliations are absolutely
continuous with respect to the invariant measure \cite[Theorem 6.2,Theorem 7.5]{Le02}
and adapting the formulation of Liverani and Wojtkowski \cite{LW95} to the
infinite measure case, he proved a local ergodicity property \cite[Theorem 8.5]{Le02}
 and as consequences its global ergodicity \cite[Theorem 8.5]{Le02}.
The definition of ergodicity used in these two results is that the Birkhoff average are constant
almost everywhere for all integrable functions, a weaker definition than the usual one,
for systems with invariant finite measure. However Lenci also proved that
$(M_3, T_3, \mu)$ is ergodic \cite[Proposition 8.11]{Le02}, which implies that
the billiard map is ergodic in the sense that invariant sets are measurably indecomposable.

\begin{remark}
It is not difficult to see that $T_3^n$ is ergodic, for all positive integer $n$.
Indeed, it is just repeated the argument in the proof of \cite[Proposition 8.11]{Le02} .
\end{remark}

Next we introduce auxiliary first return maps associated with
the billiard system we are considering. Recall that we already
defined $T_4: M_4 \rightarrow M_4$ that corresponds for the bouncing at
the vertical wall $\Le$ and $T: M \rightarrow M$ corresponding to the
bouncing at the dispersing part $\U$ of $\D$. We point out that a priori,
the bouncing at the dispersing part $\U$ contains the most of the
chaotic behavior, but this gives a system with an infinity measure,
one of the major difficulties in analyzing this billiard system.
To bypass this difficulty, we consider as well the bouncing at
the vertical wall, that gives a finite measure system. Thus
 we set $M_5 = M \cup M_4$, that consider the bouncing at the dispersing part
and the bouncing at the vertical wall and denote $T_5:M_5\to M_5$ the return map to $M_5$.
Note that by construction $M_5$ comes from bouncing at a global cross section to
the billiard, constituted by the union of $\U$ and $\Le$, the dispersing and the 
vertical wall respectively.

Thus $T_5$ is a billiard map that preserves the infinite measure
$d\mu = \cos\varphi drd\varphi$ defined at $M_5$.
When no ambiguity exists, we will use the same notation for the 
(infinite) measure invariant by the map  $T: M \rightarrow M$.

\begin{figure}[ht]
\begin{center}
\includegraphics[height=5cm, width=8cm]{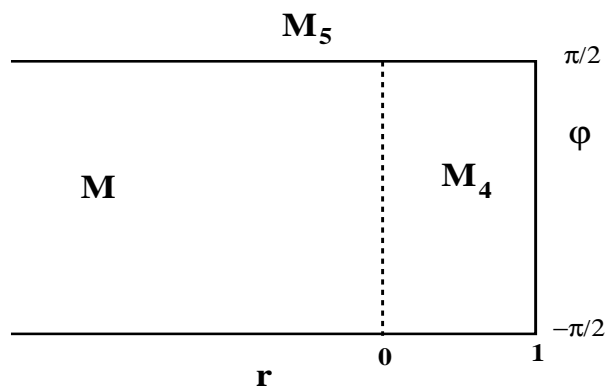}
\caption{\label{fig:m5} The phase space $M_5$. The line $r=0$,
corresponds to the vertex $V$ 
 and $r=1$ corresponds
to the vertex in the point $(0,0)$ of $\D$.}
\end{center}
\end{figure}

As in  \cite[Corollary 3.3]{Le02} we conclude that $T_5$ is conservative.
Moreover, since $T_4 : M_4 \rightarrow M_4$ is ergodic and it is induced by $T_5$,
so is $T_5$.

The map $T_5$ describes all the dynamics of our system and this is the map
we shall deal with from now on.

\section{K-automorphisms and proof of Theorem \ref{main:ksystem}} \label{sec:mainksystem}
 
\begin{definition} \label{def:ksystem}
 Let $(X,\mathcal{B}, \mu)$ be an infinite $\sigma$-finite measure space and
  $\F: X \rightarrow X$ an automorphism. We say that
$\F$ is an {\it infinite K-automorphism} if there exists a sub-$\sigma$-algebra $\mathcal{K}
\subset \mathcal{B}$ such that
\begin{itemize}
 \item [(i)] $\F\K \supset \K$;
\item [(ii)] $\bigvee_{n=0}^{\infty} \F^n\K = \B \mod \mu$;
\item [(iii)] $\bigcap_{n=0}^{\infty} \F^{-n}\K = \mathcal{N} = \{\emptyset, X\} \mod \mu$.
\end{itemize}
\end{definition}

The next proposition is an extension for K-automorphisms 
of a similar result for ergodic maps in spaces with infinite 
 measure. See \cite[p. 42]{Aa97}.

\begin{proposition} \label{prop:kinduced}
Let $\F$ be an ergodic measure preserving map of the $\sigma$-finite
measure space $(X,\B,\mu)$ and suppose there exists $E \in \B$ such
that $\mu(E) > 0$ and
$\bigcup_{n=0}^{\infty} \F^{-n}E = X \mod \mu$.
If $\F_E$ is a (finite or infinite) K-automorphism  
then $\F$ is a (finite or infinite) K-automorphism.
\end{proposition}

\begin{proof}
Since $\F_E$ is a K-automorphism, there exists $\K_E$, sub-$\sigma$-algebra of
$\B_E = \B \cap E$ such that
\begin{itemize}
 \item[(1)]  $\F_E\K_E \supset \K_E$;
\item[(2)]  $\bigvee_{n=0}^{\infty} \F^n_E\K_E = \B_E \mod \mu_E$;
\item[(3)]  $\bigcap_{n=0}^{\infty} \F^{-n}_E\K_E = \mathcal{N}_E = \{\emptyset, E=X\cap E\} \mod \mu_E$.
\end{itemize}

 Let $\K = \bigvee_{i= -\infty}^{0} \F^i\K_E$. Then

(i) \;\; $\F\K = \bigvee_{i= -\infty}^{0} \F^{i+1} \K_E =
 \bigvee_{i= -\infty}^{1} \F^i\K_E
\supset \bigvee_{i= -\infty}^{0} \F^i\K_E = \K$.

(ii) \;\; $\bigvee_{n=0}^{\infty} \F^n\K = \bigvee_{n=0}^{\infty}\F^n\bigvee_{i=-\infty}^{0}
\F^i\K_E = \bigvee_{i=-\infty}^{\infty}\F^i\K_E = \B \mod \mu$; 
since, by condition (2), $\B_E \subset \displaystyle\bigvee_{i=-\infty}^{\infty}\F^i\K_E$, 
and given $A \in \B$,
\[
 A = A \cap X = A \cap \bigcup_{i=0}^{\infty} \F^{-i}E = \bigcup_{i=0}^{\infty} (A \cap \F^{-i}E) =
 \bigcup_{i=0}^{\infty} \F^{-i}(\F^i A \cap E).
\]
Because $\F^iA\cap E \in \B_E$, it follows that
 $\F^{-i}(\F^iA \cap E) \in \F^{-i}\B_E$, so
$A \in \bigvee_{i=-\infty}^{\infty}\F^i\K_E$.
Thus $\B \subset \bigvee_{i=-\infty}^{\infty}\F^i\K_E$.

(iii) \;\; We must show that 
$\bigcap_{n=0}^{\infty}\F^{-n}\K = \bigcap_{n=0}^{\infty} \F^{-n} \bigvee_{i=-\infty}^0\F^i\K_E = \mathcal{N} \mod \mu.$
To do this, we just need to show that
$\bigcap_{n=0}^{\infty}\F^{-n}\K  \subset \mathcal{N} \mod \mu.$

Let $A \in \displaystyle\bigcap_{n=0}^{\infty} \F^{-n} \bigvee_{i=-\infty}^0\F^i\K_E$ and 
suppose that $\mu(A) > 0$. 
Furthermore, we may suppose that $\mu(A \cap E)>0$ 
because, if not, we take $A^c$ 
as the set.

Then $A \in \bigcap_{n=0}^{\infty} \F^{-n} \bigvee_{i=-\infty}^0\F^i\K_E
\Rightarrow A \cap E \in \bigcap_{n=0}^{\infty} \F^{-n} \bigvee_{i=-\infty}^0\F^i\K_E
\cap E$ which is equal to $\bigcap_{n=0}^{\infty} 
\F_E^{-n} \bigvee_{i=-\infty}^0\F_E^i\K_E$, 
by the definition of $\F_E$.

By condition (1), $\bigvee_{i=-\infty}^0\F_E^i\K_E = \K_E$ ,
so
$\bigcap_{n=0}^{\infty} \F_E^{-n} \bigvee_{i=-\infty}^0\F_E^i\K_E =
\bigcap_{n=0}^{\infty} 
\F_E^{-n} \K_E = \mathcal{N}_E \mod \mu_E$, by condition (3).
Then
$\mu_E((A\cap E)\bigtriangleup E) = 0$ and $\mu_E(A^c \cap E) = 0$.
Thus
$\mu(A^c \cap E) = 0$.

Also, $\mu_E(\F_E^k(A^c \cap E)) = \mu_E(\F_E^kA^c \cap E) = 0$, 
for all $k \geq 0$.
Then $\mu_E(\F^jA^c \cap E) = 0$, for all $j \geq 0$,
 because 
$\mu_E(\F^jA^c \cap E) \leq \sum_{k=1}^j \mu_E(\F_E^kA^c \cap E) = 0$.
Since $A^c = A^c \cap X =  A^c \cap \bigcup_{j=0}^{\infty} \F^{-j}E 
= \bigcup_{j=0}^{\infty}\F^{-j}(\F^j A^c \cap E)$, we get
 $\mu(A^c)=0$
hence $\bigcap_{n=0}^{\infty}\F^{-n}\K  \subset \mathcal{N} \mod \mu$. 
\end{proof}

We also need the following theorem due to Pesin and 
Katok and Strelcyn:

\begin{theorem} {\rm (Pesin \cite[Theorem 7.2]{Pe77} , Katok and Strelcyn \cite[Theorem 13.1]{Ka86})}  \label{th:pesin}
Let $\V$ be a finite union of compact Riemannian manifolds
 $\V_1, \V_2 \ldots \V_s$
(possibly with boundaries and corners), 
all of them with dimension $d \geq 2$, glued along finitely many
$C^1$ submanifolds of positive codimension and $\F$ a map on $\V$ 
preserving a Borel probability measure $\mu$, both satisfying the Katok
and Strelcyn conditions indicated in \cite[Section 1.1]{Ka86} .
Suppose that 
\[\Sigma(\F) = \{x \in V : \mbox{the Lyapunov exponents on $V$ are non-zero} \}\]
has positive $\mu$-measure.
Then there exist sets $\Sigma_i \subset \Sigma(\F)$, 
$i=0,1,2,\ldots$, such that 
\begin{enumerate}
\item[(1)] $\Sigma(\F) = \bigcup_{i \geq 0}\Sigma_i$, $\Sigma_i \cap \Sigma_j = \emptyset$ for
$i \neq j$, $i,j = 0,1,2,\ldots$;
\item[(2)]$\mu(\Sigma_0)=0$, $\mu(\Sigma_i) > 0$, for $i> 0$;
\item[(3)] for $i > 0$: $\F(\Sigma_i) = \Sigma_i$, $\F|\Sigma_i$
is ergodic;
\item[(4)] for $i>0$, there exists a splitting $
\Sigma_i = \bigcup_{j=1}^{n_i} \Sigma_{i}^{j}, n_i \in \mathbb{Z}^+$
such that 
\begin{enumerate}
\item[(a)]$\Sigma_i^{j_1}\cap \Sigma_i^{j_2}= \emptyset$ for $j_1 \neq j_2$;
\item[(b)] $\F(\Sigma_i^j) = \Sigma_i^{j+1}$ para $j=1,2,\ldots, n_i-1$,
 $\;\F(\Sigma_i^{n_i}) = \Sigma_i^{1}$;
\item[(c)] $\F^{n_i}|\Sigma_i^1$ is a finite K-automorphism.
\end{enumerate}
\end{enumerate}
\end{theorem}

Returning to the billiard map case:

\begin{lemma} \label{lemma:tmkaut}
 Let $M_4$ be the phase space associated to the rebounds in the vertical wall.
Then $T_{4}$ is a finite K-automorphism.
\end{lemma}
\begin{proof}

We know that $T_{4}^n$ is ergodic for all $n \geq 1$. Also, the Lyapunov
exponents for $T_4$ are non-zero. So we may apply Theorem \ref{th:pesin}.
However, by the ergodicity of $T^n_4$, all the decompositions are trivial and
 we get that $T_4$ is a K-automorphism.

\end{proof}

From Proposition \ref{prop:kinduced}, it follows that $(M_5,T_5,\mu)$ an 
infinite K-system, concluding the proof of Theorem \ref{main:ksystem}.



\section{Geometric conditions}
This section and the next one are inspired in the 
analysis for trajectories in a finite cusp studied 
by Chernov and Markarian \cite{CM07}.

We are in the same setting as in the previous sections.
Fix $N_0 >\!> 1$.
 Let us study the behavior of a trajectory  that leaving $\Le$
(with coordinates $(r,\vfi)$ in $M_4$),
enters in the cusp, and comes back after $N > N_0$ rebounds. 
In order to do this, we shall adopt a new system of coordinates from now on.
  Let $x_n \in [0, \infty)$,  $0 \leq n \leq N$, be the $x$-coordinate associated to the
   $n$-th rebound on $\U$,
(where $x_0 = 0$, leaving $\Le$),and $\ga_n \in [0,\pi/2]$, $0 \leq n \leq N$, the positive angle
between the trajectory and the tangent at the point of collision with coordinate $x_n$
($\ga_0 = \pi/2 -|\vfi|$).

Define
\[
x_{N_2} := \max{\{x_n :n=1,2,\ldots,N\}},
\]
that is, the $x$-coordinate of the most interior point
inside the cusp.

If $n \leq N_2 - 1$ then
\begin{equation} \label{eq:gamma_n+1}
\gamma_{n+1} = \gamma_n + \tan^{-1} |f^{\prime}(x_n)| +
                                                      \tan^{-1} |f^{\prime}(x_{n+1})|
 \end{equation}
\begin{equation} \label{eq:xn+1}
x_{n+1} = x_n + \frac{f(x_n)+f(x_{n+1})}{\tan(\gamma_n +
                                                     \tan^{-1} |f^{\prime}(x_n)|)}.
 \end{equation}
 If $n \geq N_2$ then
\begin{equation*}
 \gamma_{n} = \gamma_{n+1} + \tan^{-1} |f^{\prime}(x_n)| +
                                                      \tan^{-1} |f^{\prime}(x_{n+1})|
\end{equation*}
\begin{equation*}
 x_n = x_{n+1} + \frac{f(x_n)+f(x_{n+1})}{\tan(\gamma_{n+1} +
                                                  \tan^{-1} |f^{\prime}(x_{n+1})|)}.
\end{equation*}

\begin{lemma} \label{lema:N2}
Using the notation above, $|N_2 - N/2| = O(1).$
\end{lemma}

\begin{proof}
Suppose that, without lost of generality,  $x_{N_ 2+1} \geq x_{N_2-1}$. Then
\begin{equation*} \gamma_{N_2} = \gamma_{N_2-1} + \tan^{-1} |f^{\prime}(x_{N_2-1})| +
                                                \tan^{-1} |f^{\prime}(x_{N_2})|.
 \nonumber \end{equation*}
On the other hand,
\begin{equation*} \gamma_{N_2} = \gamma_{N_2+1} + \tan^{-1} |f^{\prime}(x_{N_2+1})| +
                                                \tan^{-1} |f^{\prime}(x_{N_2})|.
\nonumber \end{equation*}
So,
\begin{eqnarray*}
 \gamma_{N_2-1} + \tan^{-1} |f^{\prime}(x_{N_2-1})| & = &
             \gamma_{N_2+1} + \tan^{-1} |f^{\prime}(x_{N_2+1})| \nonumber \\
   & \leq &    \gamma_{N_2+1} + \tan^{-1} |f^{\prime}(x_{N_2-1})|, \nonumber
\end{eqnarray*}
That is,
\begin{equation*}
 \gamma_{N_2-1} \leq \gamma_{N_2+1}. \nonumber
 \end{equation*}

Now, we must show that $x_{N_2 - i} \leq x_{N_2 + i}$ and
$\gamma_{N_2-i} \leq \gamma_{N_2+i}$, for all $i= 1, 2, \ldots $ while the collisions
remain inside the cusp. Indeed, suppose that for $i$ it is true and we shall show it
for $i+1$. Then
\begin{equation*} x_{N_2-i} = x_{N_2-(i+1)} + \frac{f(x_{N_2-i})+f(x_{N_2-(i+1)})}
{\tan(\gamma_{N_2-(i+1)} + \tan^{-1} |f^{\prime}(x_{N_2-(i+1)})|)}
\nonumber \end{equation*}
\begin{equation*} x_{N_2+i} = x_{N_2+(i+1)} + \frac{f(x_{N_2+i})+f(x_{N_2+(i+1)})}
{\tan(\gamma_{N_2+(i+1)} + \tan^{-1} |f^{\prime}(x_{N_2+(i+1)})|)}
\nonumber \end{equation*}
\begin{equation*} \gamma_{N_2-i} = \gamma_{N_2-(i+1)} +
 \tan^{-1} |f^{\prime}(x_{N_2-(i+1)})| +  \tan^{-1} |f^{\prime}(x_{N_2-i})|
\nonumber \end{equation*}
\begin{equation*} \gamma_{N_2+i} = \gamma_{N_2+(i+1)} +
 \tan^{-1} |f^{\prime}(x_{N_2+(i+1)})| +  \tan^{-1} |f^{\prime}(x_{N_2+i})|.
 \nonumber \end{equation*}

By the induction hypothesis, the following holds
\begin{eqnarray*}
&&\gamma_{N_2-(i+1)} + \tan^{-1} |f^{\prime}(x_{N_2-(i+1)})| +
                                     \tan^{-1} |f^{\prime}(x_{N_2-i})| \nonumber \\
 & \leq &
\gamma_{N_2+(i+1)} + \tan^{-1} |f^{\prime}(x_{N_2+(i+1)})| +
                                   \tan^{-1} |f^{\prime}(x_{N_2+i})| \nonumber \\
& \leq & \gamma_{N_2+(i+1)} + \tan^{-1} |f^{\prime}(x_{N_2+(i+1)})| +
                                  \tan^{-1} |f^{\prime}(x_{N_2-i})|. \nonumber
\end{eqnarray*}
So
\begin{equation} \label{eq:gamma}
\gamma_{N_2-(i+1)} + \tan^{-1} |f^{\prime}(x_{N_2-(i+1)})|  \leq
\gamma_{N_2+(i+1)} + \tan^{-1} |f^{\prime}(x_{N_2+(i+1)})| .
\end{equation}

Using the hypothesis of induction and by (\ref{eq:gamma}), we get
\begin{eqnarray*}
&& x_{N_2-(i+1)} + \frac{f(x_{N_2-i})+f(x_{N_2-(i+1)})}
{\tan(\gamma_{N_2-(i+1)} + \tan^{-1} |f^{\prime}(x_{N_2-(i+1)})|)} \nonumber \\
& \leq &  x_{N_2+(i+1)} + \frac{f(x_{N_2+i})+f(x_{N_2+(i+1)})}
{\tan(\gamma_{N_2+(i+1)} + \tan^{-1} |f^{\prime}(x_{N_2+(i+1)})|)} \nonumber \\
& \leq &  x_{N_2+(i+1)} + \frac{f(x_{N_2+i})+f(x_{N_2+(i+1)})}
{\tan(\gamma_{N_2-(i+1)} + \tan^{-1} |f^{\prime}(x_{N_2-(i+1)})|)}. \nonumber
\end{eqnarray*}
Thus
\begin{equation*}
x_{N_2-(i+1)} + f(x_{N_2-(i+1)}) \leq
                                 x_{N_2+(i+1)} + f(x_{N_2+(i+1)}). \nonumber
\end{equation*}

Since $x_i \geq 1$, for all $i=1,2,\ldots, N$, 
\begin{center}
$x_{N_2-(i+1)} \leq x_{N_2+(i+1)} \;\;$ e 
$\;\;\gamma_{N_2-(i+1)} \leq \gamma_{N_2+(i+1)}$,
\end{center}
as we wish to demonstrate. Thus $|N_2 - N/2| = O(1).$
\end{proof}

Let us now split the trajectories going through the 
cusp in three regions. For this, we choose
$\bar{\gamma}$ sufficiently small, that the exact value
is not important,
 e.g. $\bar{\gamma} = 10^{-10}$.
 This choice allows us to make estimates in three different regions,
 defining
\begin{eqnarray*}
N_1 & = & \max \{ n < N_2; \gamma_n \leq \bar{\gamma}\} \nonumber \\
N_3 & = & \min \{ n > N_2; \gamma_n \leq \bar{\gamma}\}. \nonumber
\end{eqnarray*}
We call the series of rebounds between $1$ and $N_1$ the \emph{entering period},
between $N_1$ and $N_3$  the \emph{turning period} and between
$N_3$ and $N$ the \emph{exiting period}. 
Furthermore, consider $x_1$ large enough, e.g. $x_1 > 10^6$.

From now on, we use the table $\D$ 
defined by $\displaystyle f(x) = (x+1)^{-1}$.
Until the end of this section, we use the following change of variables:
\[
t_n = x_n +1, \;\; \forall 1 \leq n \leq N.
\]

\begin{lemma} \label{lema:assintotico}
 We have
\[
N_1 \asymp N_2-N_1 \asymp N_3-N_2 \asymp N-N_3 \asymp N,
\]
so all segments have size of order $N$.
Moreover
\begin{equation} \label{eq:x1}
x_1 \asymp N^{\frac{1}{6}}\;\;\; e \;\;\; x_{N_2}\asymp N^{\frac{1}{2}} 
\end{equation}
and
\begin{equation*}
x_n \asymp n^{\frac{1}{3}}N^{\frac{1}{6}} \;\;\; \forall n = 2, \ldots, N_1.
\end{equation*}
Also
\begin{equation} \label{eqn:gamma1}
\gamma_1 = O(N^{-1/3})\;\;\; e \;\;\; \gamma_2 \asymp N^{-1/3}
\end{equation}
and
\begin{equation*}
\gamma_{n} \asymp  n^{\frac{1}{3}}N^{-\frac{1}{3}} \;\;\;
                                                             \forall n = 2, \dots, N_1.
\end{equation*}
\end{lemma}

\begin{proof}

For each $n=1,2,\ldots,N_1$, we define $\omega_n =
 \dfrac{\gamma_{n}}{\frac{1}{ |f^{\prime} (t_{n})|}}$.
 Using the definition of $f$, we obtain 
$\omega_n = \gamma_{n} t_{n}^2$. Also we define
$u_n=\displaystyle \frac{t_n}{t_{n+1}}$. 
Multiplying (\ref{eq:gamma_n+1}) by $t_{n+1}^2$
and expanding $\tan^{-1}$ in its Taylor's series we get
\begin{equation} \label{eq:omega}
\omega_{n+1} = \frac{\omega_{n} + 1}{u_{n}^2} + 1 + O(x_{n}^{-4}).
\end{equation}
From (\ref{eq:gamma_n+1}), we get
\begin{equation} \label{eq:ga_1-ga_n}
\gamma_1 + \frac{1}{t_1^2} + \frac{2}{t_2^2} + \ldots + \frac{2}{t_{n-1}^2} +
\frac{1}{t_n^2} + O\left(\sum_{i=1}^n t_i^{-6}\right) = \gamma_n \leq \frac{\pi}{2}.
\end{equation}
Thus
\begin{equation} \label{eq:sum_x_i}
\sum_{i=1}^n t_i^{-2} = O(1).
\end{equation}

From equation (\ref{eq:omega}),we obtain
\begin{equation} \label{eq:winf}
\omega_n > 2n - 2.
\end{equation}
From (\ref{eq:xn+1}) and using the fact that $\tan x > x$, we have
\begin{equation}  \label{eq:supbound}
\frac{1}{u_n} < 1 + \frac{2}{\omega_n + 1} (1 + O(t_n^{-6})).
\end{equation}

Replacing (\ref{eq:supbound}) in (\ref{eq:omega}):
\begin{eqnarray} \label{eq:wn+1sup}
\omega_{n+1} & < & 1 + (\omega_n + 1)\left(1 + \frac{2}{\omega_n + 1}
           (1 + O(t_n^{-6}))\right)^2 + O(t_n^{-4}) \nonumber \\
 & = & 6 + \omega_n + \frac{4}{\omega_n + 1} + O(t_n^{-4})
 < 6 + \omega_n + \frac{4}{2n - 1} + O(t_n^{-4})\nonumber.
\end{eqnarray}
So
\begin{equation} \label{eq:wsup}
\omega_n < 6n + 2\ln n + O(1).
\end{equation}
From (\ref{eq:winf}) and (\ref{eq:wsup}) we conclude that $\omega_n = \gamma_n t_n^{2} \asymp n.$
Since $\gamma_{N_2} \approx \pi/2$, it follows that $x_{N_2}^2 \asymp N_2 \asymp N$, by Lemma
\ref{lema:N2}.

For $n=1, 2, \ldots, N_1$,
\begin{eqnarray*}
t_{n+1} & = & t_n + \dfrac{\dfrac{1}{t_n} + \dfrac{1}{t_{n+1}}}
{\gamma_n + \dfrac{1}{t_n^2} + O\left(\left(\gamma_n +
                                        \dfrac{1}{t_n^2} \right)^3\right)}
  =  t_n + \dfrac{t_n + \dfrac{t_n^2}{t_{n+1}}}{\omega_n + 1 + t_n^2
                  O\left(\left(\gamma_n + \dfrac{1}{t_n^2} \right)^3\right)}. \nonumber
\end{eqnarray*}
Dividing by $t_n$ we obtain
\begin{eqnarray*}
\dfrac{1}{u_n} & = & 1 + \dfrac{1+ u_n}{w_n + 1}
            \left(\dfrac{1}{1 + \dfrac{t_n^2}{\gamma_n t_n^2 + 1}
          O\left(\left(\gamma_n + \dfrac{1}{t_n^2} \right)^3\right)}\right). \nonumber
\end{eqnarray*}
Since the choices of $\gamma_n < 10^{-10}$ and $x_1 > 10^6$ imply that
 $O\left(\left(\gamma_n + \dfrac{1}{t_n^2} \right)^3\right) = O(\gamma_n^3)$ 
 and because 
$\displaystyle\frac{t_n^2}{\gamma_n t_n^2 + 1}O\left(\gamma_n^3\right) = 
O\left(\gamma_n^2\right)$, we obtain 
\begin{eqnarray*}
\frac{1}{u_n} & = &  1 + \frac{1+ u_n}{w_n + 1}
                                 \left(\frac{1}{1 + O(\gamma_n^2)}\right)
        = 1 + \frac{1+ u_n}{w_n + 1}(1 + O(\gamma_n^2))\nonumber,
\end{eqnarray*}
since $O(\gamma_n^2)$ is sufficiently small.
From (\ref{eq:supbound}) we get
\begin{eqnarray} \label{eq:uninfbound}
\frac{1}{u_n} & > & 1 + \left(\frac{1}{\omega_n + 1} +
          \frac{1}{\omega_n + 3 +O(t_n^{-6})}\right)(1+O(\gamma_n^2))
 >  1 + \frac{2}{\omega_n + 3} +
                               O\left(\frac{\gamma_n^2}{n}\right) \nonumber \\
& > & 1 + \frac{2}{6n + 2\ln n + O(1)} +
                               O\left(\frac{\gamma_n^2}{n}\right) \nonumber \\
& > & \exp\left(\frac{2}{6n + 2 \ln n + O(1)} -
                         \frac{4}{(6n + 2 \ln n + O(1))^2} +
                      O\left(\frac{\gamma_n^2}{n}\right)\right), \nonumber
 \end{eqnarray}
in the last inequality we use the fact that $1 + x > \exp(x - x^2)$ for small
$x$.

Multiplying from $i=1$ to $n-1$, we get
\begin{eqnarray*}
\prod_{i=1}^{n-1} u_i^{-1} & > & \exp\left(\sum \frac{2}{6i + 2 \ln i + O(1)} -
                  \sum \frac{4}{(6i + 2 \ln i + O(1))^2} +
                    O\left(\sum \frac{\gamma_i^2}{i}\right)\right) \nonumber \\
               & > & \exp(\ln n^{1/3} - C)
               =  C^{\prime} n^{1/3}, \nonumber
\end{eqnarray*}
because
\[ \sum_{i=1}^{n-1} \frac{\gamma_i^2}{i} \leq 12\sum_{i=1}^{n-1} \frac{\gamma_i^2}{\omega_i} =
12\sum_{i=1}^{n-1} \frac{\gamma_i}{x_i^2} = O(1), \]
since $\gamma_i < \pi/2$ and $\sum_{i=1}^n t_i^{-2} = O(1)$, as obtained in (\ref{eq:sum_x_i}).
We obtain
\begin{equation} \label{eq:xnx1g}
\frac{t_n}{t_1} > C^{\prime}n^{1/3}.
\end{equation}
However
\begin{eqnarray*}
\omega_{n+1} & = & 1 + \frac{\omega_n +1}{u_n^2} + O(t_n^{-4})  \nonumber \\
          & > & 1 + (\omega_n +1)
\left(1 + \frac{2}{\omega_n +3} + O\left(\frac{\gamma_n^2}{n}\right)\right)^2
                                                   + O(t_n^{-4}) \nonumber \\
& > & 2 + \omega_n + 4\left(1 - \frac{2}{\omega_n + 3}\right) +
      4\left(\frac{1}{\omega_n + 3} - \frac{2}{(\omega_n + 3)^2}\right) +
                                     O(\gamma_n^2) + O(t_n^{-4}) \nonumber \\
& > & \omega_n + 6 -\frac{4}{2n+1} - \frac{4}{(2n + 1)^2} +
                                     O(\gamma_n^2) + O(t_n^{-4}). \nonumber
\end{eqnarray*}
This implies that
\begin{equation*}
\omega_n > 6n - 2 \ln n + O(1) \nonumber.
\end{equation*}
Thus
\begin{eqnarray*}
\dfrac{1}{u_n^2} & = &
                 \dfrac{\omega_{n+1} -1 + O(t_n^{-4})}{\omega_n +1}
  <  \dfrac{\omega_n + 1 + 4 + \dfrac{4}{\omega_n + 1} + O(t_n^{-4})}
                                                    {\omega_n + 1} \nonumber \\
 & = & 1 + \dfrac{4}{\omega_n +1} + \dfrac{4}{(\omega_n + 1)^2} +
                                  O\left(\dfrac{t_n^{-4}}{n}\right) \nonumber \\
 & < & 1 + \dfrac{4}{6n - 2 \ln n + O(1)} + \dfrac{4}{(6n - 2 \ln n + O(1))^2} +
                                  O\left(\dfrac{t_n^{-4}}{n}\right). \nonumber
\end{eqnarray*}

Multiplying from $i=1$ to $n-1$, 
\begin{eqnarray*}
\prod_{i=1}^{n-1} u_i^{-2} & < & \exp \left( \sum \frac{4}{6i - 2 \ln i + O(1)} +
                      \sum \frac{4}{(6i - 2 \ln i + O(1))^2} +
                       O\left(\sum \frac{x_i^{-4}}{i}\right)\right) \nonumber \\
               & < & \exp(\ln n^{2/3} + C)
               =  C^{\prime}n^{2/3}.    \nonumber
\end{eqnarray*}
And we obtain
\begin{equation} \label{eq:xnx1l}
\left(\frac{t_n}{t_1}\right)^2 < C^{\prime}n^{2/3}.
\end{equation}

From (\ref{eq:xnx1g}) and (\ref{eq:xnx1l}) we get
\begin{equation*}
  \frac{t_n}{t_1} \asymp n^{1/3} \nonumber.
\end{equation*}
So
\[ n \asymp \gamma_n t_n^2 \asymp \gamma_n n^{2/3}t_1^2.\]
However $\gamma_{N_1} \approx \bar{\gamma} =$ constant, then
\[ N_1 \asymp \gamma_{N_1} N_1^{2/3}x_1^2 \Rightarrow x_1 \asymp N_1^{1/6}.\]
And thus
\[ t_n \asymp n^{1/3}N_1^{1/6} \;\;
\mbox{ and } \;\;
 \gamma_n \asymp n^{1/3}N_1^{-1/3}, \]
 for all $n=2, \ldots, N_1$.

To show that $N_1 \asymp N$, just notice that at the turning period,
 i.e., $N_1 \leq n \leq N_2$, the angle $\gamma_n$ increases
from $\bar{\gamma}$ to approximately $\pi/2$ and
\[ \frac{1}{t_n^2} = \frac{\gamma_n}{\omega_n} > \frac{\bar{\gamma}}{6n+2\ln n +C}.
\]
It follows from (\ref{eq:ga_1-ga_n}) and (\ref{eq:sum_x_i}) that
\[
\sum_{n=N_1}^{N_2} (\gamma_n-\gamma_{n-1}) \geq
\sum_{n=N_1}^{N_2} \frac{C^\prime}{6n+2\ln n +C} \geq C^{\prime\prime}\ln \frac{N_2}{N_1},
\]
for some constants $C^{\prime}, C^{\prime\prime} > 0$. This implies that
$N_1 < N_2 < C^{\prime\prime\prime}N_1$, for some $C^{\prime\prime\prime}>0$.
\end{proof}

In the proof of Lemma \ref{lema:assintotico}, we obtained the following
\[
 \omega_n + 1 = 6 + \omega_n + O\left(\frac{1}{n} + \gamma^{2}_n +
t_n^{-4}\right),
\]
\[
\omega_n = 6n + O(\ln n),
\]
\[
\frac{1}{u_n} = 1 + \frac{1}{3n} + O\left(\frac{\ln n}{n^2} +
\frac{\gamma_n^2}{n} + \frac{t_n^{-4}}{n}\right),
\]
\[
u_n = 1 - \frac{1}{3n} + O\left(\frac{\ln n}{n^2} +
\frac{\gamma_n^2}{n} + \frac{t_n^{-4}}{n}\right),
\]
and we shall use the values from now on.

For $1 \leq n \leq N_2$,
let $\tau_n$ be the time between two consecutive collisions
in the billiard table:
\[ \tau_n = \frac{f(t_n) + f(t_{n+1})}{\sin(\gamma_n +
                                            \tan^{-1}(|f^{\prime}(t_n)))}. \]

Then, using the values above,
\begin{eqnarray*}
\tau_n &=& \dfrac{\dfrac{1}{t_n}+\dfrac{1}{t_{n+1}}}
           {\left(\gamma_n+\dfrac{1}{t_n^2}\right)+
        O\left(\left(\gamma_n+\dfrac{1}{t_n^2}\right)^3\right)}
       = \dfrac{t_n + \dfrac{t_n^2}{t_{n+1}}}
                            {\omega_n + 1 +O(t_n^2\gamma_n^3)} \nonumber \\
       &=& \dfrac{t_n(1+u_n)}
                     {\omega_n(1+\omega_n^{-1}+O(\gamma_n^2))}
       = \dfrac{t_n}{\omega_n}
                \frac{2+ O(n^{-1})}{1+O(n^{-1})+O(\gamma^2_n)}  \nonumber \\
       &=& \displaystyle\frac{2t_n}{\omega_n}
      \dfrac{1}{\dfrac{1+O(n^{-1})+O(\gamma^2_n)}{1+ O(n^{-1})}}
       = \displaystyle\frac{2t_n}{\omega_n}(1+O(n^{-1})+O(\gamma_n^2))     \nonumber \\
   &\asymp& n^{-2/3}N^{1/6},
\end{eqnarray*}
for $1 \leq n \leq N_2$.

Furthermore, if we denote by $K_n$ the curvature of the dispersing part
of the table at the point of collision $(r_n, \varphi_n)$,
we have that, as w enter in the cusp ($n= 1,2,\ldots, N_2$)
\begin{equation*}
K_n = \dfrac{f^{\prime\prime}(t_n)}{(1+(f^{\prime}(t_n))^2)^{3/2}}
=\dfrac{\dfrac{2}{t_n^3}}
                 {\left(1+ \left(\dfrac{1}{t_n^2}\right)^2\right)^{3/2}}
= \dfrac{\dfrac{2}{t_n^3}}
                 {\left(\dfrac{t_n^4 +1}{t_n^4}\right)^{3/2}}
= \dfrac{2}{t_n^3},
\end{equation*}
since $x_1 > 10^6$.

We also have that
\begin{eqnarray} \label{eq:tauK}
\dfrac{\tau_n K_n}{\sin \gamma_n} &=& \dfrac{2t_n\omega_n\sp{-1}
(1+O(n\sp{-1})+O(\gamma_n^2))}{\gamma_n + O(\gamma_n\sp{3})}
 \dfrac{2}{t_n^3}
= \dfrac{4}{t_n^2}\dfrac{(1+O(1/n)+O(\gamma_n^2))}
{\omega_n \ga_n(1+O(\ga_n^2))} \nonumber \\
&=&\dfrac{4}{t_n^2}(1+O(1/n)+O(\gamma_n^2))(1+O(\ga_n^2))
    = \dfrac{4}{\omega_n^2}(1 + O(n^{-1}) + O(\gamma_n^2)) \nonumber \\
&=& \dfrac{4}{(6n+O(\ln n))^2}(1 + O(n^{-1}) + O(\gamma_n^2)) \nonumber \\
&=&\dfrac{4}{36n^2 +O(n\ln n) +O((\ln n)^2)}(1 + O(n^{-1}) + O(\gamma_n^2)) \nonumber \\
&=&\dfrac{4}{36n^2(1+O(\ln n/n)+O((\ln n/n)^2))}
                                             (1 + O(n^{-1}) + O(\gamma_n^2)) \nonumber \\
&=& \dfrac{1}{9n^2}\left(1+O\left(\dfrac{\ln n}{n}\right)\right)
                                            (1 + O(n^{-1}) + O(\gamma_n^2)) \nonumber \\
&=& \dfrac{1}{9n^{2}} + O\left(\dfrac{\ln n}{n^{3}} +
 + \dfrac{\gamma_n^{2}}{n^{2}} \right).
\end{eqnarray}

Moreover
\[
\frac{\tau_{n+1}}{\tau_n} =
\frac{f(t_{n+1}) + f(t_{n+2})}{f(t_n) + f(t_{n+1})} 
\frac
{\sin(\gamma_{n} + \tan^{-1}(|f^{\prime}(t_{n})|))}
{\sin(\gamma_{n+1} + \tan^{-1}(|f^{\prime}(t_{n+1})|))} = F_1  F_2.
\]

To obtain $F_1$, we notice that
\begin{eqnarray*}
F_1 &=& \frac{\displaystyle{\frac{1}{t_{n+1}} + \frac{1}{t_{n+2}}}}
{\displaystyle{\frac{1}{t_{n+1}} + \frac{1}{t_n}}}
= \frac{\displaystyle{\frac{1}{t_{n+1}}
\left(1 + \frac{t_{n+1}}{x_{n+2}}\right)}}
{\displaystyle{\frac{1}{x_{n}}\left(1 + \frac{t_n}{t_{n+1}}\right)}} \nonumber  \\
 &=& u_n  \frac{(1+u_{n+1})}{(1+u_n)}.
\end{eqnarray*}

This last one can be computed as 
\begin{eqnarray*}
\frac{1+u_{n+1}}{1+u_n} &=&
 \frac{\displaystyle{1+1-\frac{1}{3(n+1)}+
O\left(\frac{\ln(n+1)}{(n+1)^2}+
 \frac{\ga_{n+1}^2}{n+1}+\frac{t_{n+1}^{-4}}{n+1} \right)}}
{\displaystyle{1+1-\frac{1}{3n}+
O\left(\frac{\ln n}{n^2}+
 \frac{\ga_n^2}{n}+\frac{t_n^{-4}}{n} \right)}} \nonumber \\
&=& \frac{6(n+1) - 1 + O\left(\displaystyle{\frac{\ln(n+1)}{n+1} + \ga_{n+1}^2+t_{n+1}^{-4}}\right)}
{6n - 1 + O\left(\displaystyle{\frac{\ln n}{n} + \ga_n^2+t_n^{-4}}\right)}
\frac{3n}{3(n+1)} \nonumber \\
&=&\left(1-\frac{1}{6(n+1)}+O\left(\frac{\ln n+1}{{n+1}^2} +
\frac{\ga_{n+1}^2}{n+1}+\frac{t_{n+1}^{-4}}{n+1}\right)\right)\times  \nonumber \\
&\times&\left(
1+\frac{1}{6n}+O\left(\frac{\ln n}{n^2} +\frac{\ga_n^2}{n}+\frac{t_n^{-4}}{n}\right)\right)  \nonumber \\
&=&1+\frac{5}{36n^2}\left(1+O\left(\frac{1}{n}\right)\right)+
O\left(\frac{\ln n}{n^2} +\frac{\ga_n^2}{n}+\frac{t_n^{-4}}{n}\right) \nonumber \\
&=&1+\frac{5}{36n^2}+
O\left(\frac{\ln n}{n^2} +\frac{\ga_n^2}{n}+\frac{t_n^{-4}}{n}\right).
\end{eqnarray*}

And to obtain $F_2$,
\begin{eqnarray*}
F_2 &=&\frac{\displaystyle{\ga_n+\frac{1}{t_n^2}+O(\ga_n^3)}}
{\displaystyle{\ga_{n+1}+\frac{1}{t_{n+1}^2}+O(\ga_{n+1}^3)}}
=\frac{t_n^2}{t_n^2}\frac{t_{n+1}^2}{t_{n+1}^2}
\frac{\displaystyle{\ga_n+\frac{1}{t_n^2}+O(\ga_n^3)}}
{\displaystyle{\ga_{n+1}+\frac{1}{t_{n+1}^2}+O(\ga_{n+1}^3)}}\nonumber \\
&=&\frac{t_{n+1}^2}{t_n^2}
\frac{\omega_n + 1 + O(\ga_n^2\omega_n)}{\omega_{n+1} + 1
+O(\ga_{n+1}^2\omega_{n+1})}
= \frac{1}{u_n^2} \frac{\omega_n + 1 +O(\gamma_n^2\omega_n)}
{\omega_{n+1} + 1 + O(\gamma_{n+1}^2\omega_{n+1})}   \nonumber \\
  &=& \frac{1}{u_n^2}   \frac{\omega_n + 1 +O(\gamma_n^2\omega_n)}
{\omega_n + 7 + O(\gamma_{n+1}^2\omega_{n+1} + n^{-1} + \gamma_n^2 +
                                             t_n^{-4})} \nonumber \\
&=&\frac{1}{u_n^2}\left(1-\frac{1}{n}+\frac{7}{6n^2}+
O\left(\frac{\ln n}{n^2} + \frac{\ga_n^2}{n}+\frac{t_n^{-4}}{n}\right)\right).
\end{eqnarray*}
Hence
\begin{eqnarray} \label{eq:taun+1taun}
\frac{\tau_{n+1}}{\tau_n} &=& u_n\left(1+\frac{5}{36n^2}+
O\left(\frac{\ln n}{n^2} +\frac{\ga_n^2}{n}+\frac{t_n^{-4}}{n}\right)\right)\times \nonumber \\
&\times&\frac{1}{u_n^2}\left(1-\frac{1}{n}+\frac{7}{6n^2}+
O\left(\frac{\ln n}{n^2} + \frac{\ga_n^2}{n}+\frac{t_n^{-4}}{n}\right)\right) \nonumber \\
&=&\frac{1}{u_n}\left(1+\frac{5}{36n^2}+
O\left(\frac{\ln n}{n^2} +\frac{\ga_n^2}{n}+\frac{t_n^{-4}}{n}\right)\right)\times \nonumber \\
&\times&\left(1-\frac{1}{n}+\frac{7}{6n^2}+
O\left(\frac{\ln n}{n^2} + \frac{\ga_n^2}{n}+\frac{t_n^{-4}}{n}\right)\right) \nonumber \\
&=&\left(1+\frac{1}{3n}\right)\left(1+\frac{5}{36n^2}\right)
\left(1-\frac{1}{n}+\frac{7}{6n^2}\right)
+O\left(\frac{\ln n}{n^2} + \frac{\ga_n^2}{n}+\frac{t_n^{-4}}{n}\right) \nonumber \\
&=& 1 - \frac{2}{3n} +
O\left(\frac{\ln n}{n^2} + \frac{\gamma_n^2}{n} +\frac{t_n^{-4}}{n}\right).
\end{eqnarray}

\begin{remark}
Due to the reversibility property of the billiard map, all the formulas obtained
above hold for the exiting period as well. So
\begin{equation} \label{eqn:gamman}
x_{N} \asymp N^{1/6}\;\;\;\mbox{e}\;\;\;\gamma_N = O(N^{-1/3}).
\end{equation}
During the exiting period we can use the countdown index $m=N+1-n$
obtaining asymptotic rates for
$m =N_3 -1, \ldots, N_1$, as for example, 
$x_m \asymp m^{1/3}N^{1/6}$, $\tau_m \asymp m^{-2/3}N^{1/6}$, etc.
\end{remark}

\section{Hyperbolicity} \label{sec:hyprate}

We use in this section the \emph{p-norm}, defined by
\[
\|dx\|_p = \cos \varphi |dr|,
\]
for vectors $dx \in T_xM$ of a point $x = (r,\varphi)$. 
For billiard maps, the expansion rate of unstable vectors
(i.e., in an unstable cone) in the p-norm is given by
\[
\displaystyle\frac{\|D_x\T^{n+1}(dx)\|_p}{\|dx\|_p}=\prod_{i=0}^{n}|1+\tau_iB_i|
\]
(see \cite[p.58]{CM06}). Here $B_i$ denotes
the curvature of a small arc transverse to the wave
front. For further details we suggest Chernov and Markarian's book
\cite[Chapter IV]{CM03}.
Moreover, for semi-dispersing billiards, unstable vectors are expanded monotonically
in the p-norm, and this is not necessarily true in the Euclidean norm
 (see \cite[Section 4.4]{CM06}).

 The values $B_i^+$ can be calculated inductively as
\[
B^{+}_{n+1} = \frac{2K_{n+1}}{\sin \gamma_{n+1}} +
\frac{B_{n}^{+}}{1 + \tau_n B^{+}_n}.
\]
From the equations (\ref{eqn:gamma1}) and (\ref{eqn:gamman}),
 we know that $\ga_1 =O(N^{-1/3})$ and $\ga_N
=O(N^{-1/3})$, hence they can be arbitrarily close to zero, 
which implies that the expansion rate would be extremely high.
However $B_{n+1}^+$ is an increasing function of  $B_{n}^+$ and $\dfrac{1}{\sin \gamma_{n+1}}$.
So if $\gamma_n$ increases,
$B_n^+$ decreases. In this way, we can obtain an upper bound for the expansion rate
 taking lower bounds for $\gamma_1$ and $\gamma_N$.
Thus we introduce the following assumption
\begin{equation} \label{eq:gamma1N}
\ga_1 \asymp N^{-1/3} \;\;\;\mbox{and}\;\;\; \ga_N \asymp N^{-1/3}.
\end{equation}

Let
\[
E_N= \{x \in M_4 \;\;|\;\; R(x)= N+1\},
\]
 $N > N_0$, where 
\[
R(x) = \inf\{n \geq 1 : T_ 5^nx \in M_4\},
\] 
i.e., $R(x)-1$ indicates the amount of rebounds
in the dispersing part before returning to the vertical wall
$\mathcal{L}$. So, $E_N$ is the subset of $M_4$ 
that return for the first time to $M_4$ after $N+1$ iterations of $T_5$.

The main goal of this section is to prove the following theorem

\begin{theorem} \label{teo:derivada}
For all $x \in E_N$, satisfying $\ga_1\asymp N^{-1/3}$ e $\ga_N \asymp N^{-1/3}$,
\[
 \frac{\|D_xT_5^{N+1}(dx)\|_p}{\|dx\|_p} \asymp N.
\]
\end{theorem}

Let us denote $\tau_iB_i$ by $\lambda_i$.
For $n \geq 1$,
\begin{equation} \label{eq:lambdan+1}
\lambda_{n+1} = \frac{2\tau_{n+1}K_{n+1}}{\sin\gamma_{n+1}}
+\frac{\tau_{n+1}}{\tau_n}\cdot \frac{\lambda_n}{1+\lambda_n}.
\end{equation}

\begin{lemma} \label{lemma:lambda} We have that 
\begin{eqnarray*}
\lambda_n \asymp \frac{1}{n} &,& 1 \leq n \leq N_1,    \nonumber \\
 \lambda_n \asymp \frac{1}{n} \asymp \frac{1}{N} &,& N_1 \leq n \leq N_3, \nonumber \\
\lambda_n \asymp \frac{1}{(N-n)} &,& N_3 \leq n < N. \nonumber
\end{eqnarray*}
\end{lemma}

\begin{proof}
For $1 \leq n \leq N_1$, 
$
\lambda_{n+1} > \dfrac{a}{n^2}
+\left(1- \frac{b}{n}\right) \dfrac{\lambda_n}{1+\lambda_n},
$
for some $a,b >0$. Suppose that $\lambda_n > c/n$. Then
\begin{eqnarray*}
\lambda_{n+1} &>& \frac{a}{n^2}
      +\left(1- \frac{b}{n}\right) \frac{c/n}{1+c/n}
    = \frac{a}{n^2} +\left(1- \frac{b}{n}\right)
                               \left(\frac{c}{n+c}\right) \nonumber \\
   &=& \frac{a}{n^2} + \frac{c}{n+c} - \frac{bc}{(n+c)n}
   = \frac{a(n+c) + cn^2-bcn}{(n+c)n^2}             \nonumber     \\
    &=& \frac{c + \left(a-bc + ac/n\right)/n}{n+c}.
\end{eqnarray*}

 If $c>0$ is small enough, the expression in parenthesis is positive and
$
\lambda_{n+1} > \dfrac{c}{n+c} > \dfrac{c}{n+1}.
$
Similarly
$
\lambda_{n+1} < \dfrac{A}{n^2}
+\left(1- \frac{B}{n}\right) \dfrac{\lambda_n}{1+\lambda_n}.
$

Supposing $\lambda_n < C/n$, we get
$
\lambda_{n+1} < \dfrac{C + \left(A-BC + AC/n\right)/n}{n+C}.
$

If $C>0$ is large enough, the expression in parenthesis is negative for $N$ large and
$
\lambda_{n+1} < \dfrac{C}{n+C} < \dfrac{C}{n+1},
$
completing the induction.

For $N_1 \leq n \leq N_3$, 
$
\lambda_{N_1} \asymp \dfrac{1}{N}$ e
$
\tau_n \asymp n^{-2/3}N^{1/6} \asymp N^{-1/2}.
$
So
$
B^{+}_{N_1} = \dfrac{\lambda_{N_1}}{\tau_{N_1}} \asymp N^{-1/2}.
$
We have that
\[
K_{n+1} \asymp N^{-3/2} \Rightarrow \exists \; a,A>0 \;\; \mbox{such that}
\;\; aN^{-3/2} \leq K_{n+1} \leq AN^{-3/2},
\]
\[
\tau_n \asymp N^{-1/2} \Rightarrow \exists \; b,B>0 \;\; \mbox{such that}
\;\; bN^{-1/2} \leq \tau_n \leq BN^{-1/2},
\]
\[
B_{N_1}^{+} \asymp N^{-1/2} \Rightarrow \exists\; c,C>0 \;\;\mbox{such that}
\;\; cN^{-1/2} \leq B_{N_1}^{+} \leq CN^{-1/2}.
\]
Moreover
\[ 2 \leq \frac{2}{\sin \gamma_{n+1}} \leq \frac{2}{\sin \bar{\gamma}}
=: G, \;\;\; \forall N_1 \leq n \leq N_3.
\]
So
\[
B_{N_1+1}^{+} = \frac{2K_{N_1+1}}{\sin\gamma_{N_1+1}}+\frac{B^{+}_{N_1}}
                                   {1+\tau_{N_1}B^{+}_{N_1}}
    \leq GAN^{-3/2} + CN^{-1/2}.
\]
\begin{eqnarray*}
B_{N_1+2}^{+} &=& \frac{2K_{N_1+2}}{\sin\gamma_{N_1+2}}+\frac{B^{+}_{N_1+1}}
                              {1+\tau_{N_1+1}B^{+}_{N_1+1}}
   \leq GAN^{-3/2} + B_{N_1+1}^{+}  \nonumber  \\
   & \leq & 2GAN^{-3/2} + CN^{-1/2}.\nonumber
\end{eqnarray*}
Thus
\begin{eqnarray*}
B_{n}^{+} &=&  \frac{2K_{n}}{\sin\gamma_{n}}+\frac{B^{+}_{n-1}}
                              {1+\tau_{n-1}B^{+}_{n-1}}
 \leq  (n- N_1)GAN^{-3/2} + CN^{-1/2}    \nonumber  \\
& \leq & DNGAN^{-3/2} + CN^{-1/2}
 =  (DGA+C)N^{-1/2} \nonumber  \\
& =& EN^{-1/2}.\nonumber
\end{eqnarray*}

On the other hand
\[
B_{N_1+1}^{+} = \frac{2K_{N_1+1}}{\sin\gamma_{N_1+1}}+\frac{B^{+}_{N_1}}
                                   {1+\tau_{N_1}B^{+}_{N_1}}
 \geq  2aN^{-3/2} + \frac{cN^{-1/2}}{1+BEN^{-1}}.
\]
\begin{eqnarray*}
B_{N_1+2}^{+} &=& \frac{2K_{N_1+2}}{\sin\gamma_{N_1+2}}+
   \frac{B^{+}_{N_1+1}}{1+\tau_{N_1+1}B^{+}_{N_1+1}} \nonumber  \\
& \geq & 2aN^{-3/2} + \frac{1}{(1+BEN^{-1})}
\left(2aN^{-3/2} + \frac{cN^{-1/2}}{1+BEN^{-1}}\right) \nonumber \\
&=& 2aN^{-3/2}\left(1 + \frac{1}{(1+BEN^{-1})}\right) +
        \frac{cN^{-1/2}}{(1+BEN^{-1})^2}. \nonumber
\end{eqnarray*}
\begin{eqnarray*}
B_n^{+} & =&  \frac{2K_{n}}{\sin\gamma_{n}}+
   \frac{B^{+}_{n-1}}{1+\tau_{n-1}B^{+}_{n-1}} \nonumber \\
& \geq&  2aN^{-3/2}\left(\sum_{i=0}^{n-N_1-1}
\frac{1}{(1+BEN^{-1})^i}\right) + \frac{cN^{-1/2}}{(1+BEN^{-1})^{n-N_1}}.\nonumber
\end{eqnarray*}

There exist constants $f > 0$ and $h > 0$ such that $\sum_{i=0}^{n-N_1-1}
\frac{1}{(1+BEN^{-1})^i} \geq fN\;$ 
and also $\;\frac{1}{(1+BEN^{-1})^{n-N_1}}
\geq h$, because $\;\sum_{i=0}^{n-N_1-1}
\frac{1}{(1+BEN^{-1})^i} \geq \sum_{i=0}^{n-N_1-1}
\frac{1}{(1+BEN^{-1})^{N_3-N_1}} \asymp N\;$, and
$\;\frac{1}{(1+BEN^{-1})^{n-N_1}}$ is a bounded sequence.

So
\begin{eqnarray*}
B_{n}^{+} &\geq& 2afN^{-1/2} + chN^{-1/2} \nonumber \\
 &=& (2af+ch)N^{-1/2} \nonumber  \\
 & = & eN^{-1/2}.\nonumber
\end{eqnarray*}
Thus $B^{+}_n \asymp N^{-1/2}$,
and therefore $\lambda_n = B^{+}_n\tau_n \asymp N^{-1/2}N^{-1/2} = N^{-1}$.

For $ N_3 \leq n < N$, using the reversibility property of the billiard map,
\[
\lambda_{m-1} = \frac{2\tau_{m-1}K_{m-1}}{\sin\gamma_{m-1}}
+\frac{\tau_{m-1}}{\tau_m}\cdot \frac{\lambda_m}{1+\lambda_m},
\]
for $m= N+1 -n$. In particular,
\[
\frac{a}{m^2} < \frac{2\tau_{m-1}K_{m-1}}{\sin\gamma_{m-1}}
< \frac{A}{M^2}
\;\;\mbox{ e }\;\;
1 + \frac{b}{m} < \frac{\tau_{m-1}}{\tau_m} < 1 + \frac{B}{m},
\]
for some $0<a<A<\infty$ and $0<b<B<\infty$.

Supposing $\lambda_m > c/m$,
\begin{eqnarray*}
\lambda_{m-1} &>& \frac{a}{m^2} + \left(1+ \frac{b}{m}\right)
\frac{c/m}{1 + c/m}   \nonumber  \\
&=& \frac{c + [a + bc - c- c^2 + (ac -a - bc -ac/m)/m]/(m+c)}{m-1} \nonumber
\end{eqnarray*}
If $c>0$ is small enough, the expression between brackets is positive, for $m$ large, 
and we obtain that $\lambda_{m-1} > c/(m-1)$.

Supposing that $\lambda_m < C/m$, 
\begin{eqnarray*}
\lambda_{m-1} &<& \frac{A}{m^2} + \left(1+ \frac{B}{m}\right)
\frac{C/m}{1 + C/m}  \nonumber   \\
&=& \frac{C + [A + BC - C- C^2 + (AC -A - BC -AC/m)/m]/(m+C)}{m-1} \nonumber
\end{eqnarray*}
If $C>0$ is large enough, the expression between brackets is negative, for $m$ large,
and we obtain $\lambda_{m-1} < C/(m-1)$,
completing the proof.
\end{proof}

Lemma \ref{lemma:lambda} implies that $\displaystyle\sum_{n=1}^{N-1}\lambda_n^2=O(1)$.
Therefore, for $1 \leq N^{\prime} < N^{\prime\prime} \leq N$, 
\begin{equation} \label{eq:prodsum}
\prod_{n=N^{\prime}}^{N^{\prime\prime}-1}(1+\lambda_n) = 
\exp\left(\sum_{n=N^{\prime}}^{N^{\prime\prime}-1}\ln(1+\lambda_n)\right)
\asymp
\exp\left(\sum_{n=N^{\prime}}^{N^{\prime\prime}-1}\lambda_n\right).
\end{equation}
in the turning period, we have that $\sum_{n=N_1}^{N_3-1}\lambda_n \asymp 1$,
showing that the expansion during this period is negligible.

\begin{lemma} \label{lema:lambda_1_N1}
For all $x \in E_N$ satisfying (\ref{eq:gamma1N}),
$\prod_{n=1}^{N_1}(1+\la_n) \asymp N^{2/3}$.
\end{lemma}

\begin{proof}
According to the equation (\ref{eq:prodsum}), 
 it is sufficient to show that
\[
\la_n = \frac{2}{3n} + \chi_n; \;\;\mbox{where}\;\;
\sum_{n=1}^{N_1}\chi_n = O(1).
\]

We have that, by (\ref{eq:lambdan+1})
\[
\la_{n+1} = \frac{2}{9n^2} + a_n + \left(1 - \frac{2}{3n} + b_n\right)
\frac{\la_n}{1+\la_n},
\]
where
\[
a_n = O\left(\frac{\ln n}{n^3} + \frac{\ga_n^2}{n^2}\right) \;\;\mbox{e}
\;\; b_n = O\left(\frac{\ln n}{n^2}+\frac{\ga_n^2}{n}+\frac{x_n^{-4}}{n}
\right),
\]
are relative to the equations (\ref{eq:tauK}) and (\ref{eq:taun+1taun}).

Note that $|a_n| \leq c/n^2$ and $|b_n| \leq c/n$, for some $c > 0$
small enough.

Take
\[
\la_n = 2 \frac{1+Z_n}{3n}.
\]
We get that
\begin{eqnarray*}
2\frac{1+Z_{n+1}}{3(n+1)} &=& \frac{2}{9n^2} + a_n +
                           \left(1 - \frac{2}{3n} + b_n\right) \times \nonumber \\
&& \times
                            \left(\frac{2}{3n} + \frac{2Z_n}{3n}\right)
     \left(1 - \frac{2}{3n} - \frac{2Z_n}{3n} + O\left(\frac{1}{n^2}
                      +\frac{Z_n^2}{n^2}\right)\right) \nonumber \\
&=&  \frac{2}{9n^2} + a_n + X_1\cdot X_2 \cdot X_3,
\end{eqnarray*}
\begin{eqnarray*}
X_2\cdot X_3 &=& \frac{2}{3n} + \frac{2Z_n}{3n} - \frac{4}{9n^2} -
\frac{8Z_n}{9n^2} - \nonumber \\
&& -
\frac{4Z_n^2}{9n^2} + O\left(\frac{1}{n^3} +
\frac{Z_n}{n^3} +\frac{Z_n^2}{n^3}+ \frac{Z_n^3}{n^3}\right) \nonumber
\end{eqnarray*}
\begin{eqnarray*}
X_1\cdot X_2 \cdot X_3 &=& \frac{2}{3n} - \frac{8}{9n^2} + \frac{2b_n}{3n}
+ \frac{2Z_n}{3n} - \frac{12Z_n}{9n^2} + \nonumber \\
&& + \frac{2b_nZ_n}{3n} -
\frac{4Z_n^2}{9n^2} +  O\left(\frac{1}{n^3} +
\frac{Z_n}{n^3} +\frac{Z_n^2}{n^3}+ \frac{Z_n^3}{n^3}\right),\nonumber
\end{eqnarray*}
Therefore
\begin{eqnarray*}
Z_{n+1} &=& R_n + Z_n \times \nonumber \\
&&
\left(1-\frac{1}{n} + b_n + O\left(\frac{1}{n^2}\right)- Z_n
\left(\frac{2}{3n}+ O\left(\frac{1}{n^2}\right)\right) +
O\left(\frac{Z_n^2}{n^2}\right)\right),\nonumber
\end{eqnarray*}
where
\[
R_n =  \frac{3}{2}na_n + b_n +O\left(\frac{1}{n^2}\right).
\]

If we fix a small $\de>0$, then for $n$ large enough
\[
|Z_{n+1}| \leq |R_n| + |Z_n|\left(1-\frac{\delta}{n}\right).
\]
Without affecting the asymptotic behavior of $Z_n$, we can assume that the upper bound 
 holds for all $n$. Using it recurrently we get that
\begin{eqnarray*}
|Z_n| &\leq& |R_n| + \sum_{k=1}^{n-1}|R_k|\prod_{i=k}^{n-1}
                            \left(1-\frac{\de}{i+1}\right) \nonumber \\
&\leq & \const\sum_{k=1}^{n}\left(|R_k|\exp\left(-\sum_{i=k}^{n}
                         \frac{\de}{(i+1)}\right)\right)  \nonumber \\
&\leq & \const\sum_{k=1}^{n}|R_k|(k/n)^{\de}.
\end{eqnarray*}
Then
\begin{eqnarray*}
\sum_{n=1}^{N_1}|\chi_n| &\leq& \sum_{n=1}^{N_1}|Z_n|/n \nonumber \\
    &\leq&        \const\sum_{n=1}^{N_1}
                        \sum_{k=1}^{n}|R_k|k^{\de}/n^{\de+1} \nonumber \\
&\leq&         \const\sum_{k=1}^{N_1}
                     |R_k| \sum_{n=k}^{N_1}k^{\de}/n^{\de+1}  \nonumber \\
&\leq&           \const\sum_{k=1}^{N_1}|R_k|.
\end{eqnarray*}

The last sum is uniformly bounded on $N$, which completes the proof.
\end{proof}

\begin{lemma} \label{lema:lambda_N3_N}
For all $x \in E_N$ satisfying (\ref{eq:gamma1N})
$\prod_{n=N_3}^{N}(1+\la_n) \asymp N^{1/3}$.
\end{lemma}

\begin{proof}
 It is sufficient to show that, for $m=N-n+1$,
\[
\la_m = \frac{1}{3m} + \chi_m; \;\;\mbox{where}\;\;
\sum_{m=2}^{N-N_3}\chi_m = O(1).
\]

We have that
\[
\la_{m-1} = \frac{2}{9m^2} + a_m + \left(1 + \frac{2}{3m} + b_m\right)
\frac{\la_m}{1+\la_m},
\]
where
\[
a_m = O\left(\frac{\ln m}{m^3} + \frac{\ga_m^2}{m^2}\right) \;\;\mbox{e}
\;\; b_m = O\left(\frac{\ln m}{m^2}+\frac{\ga_m^2}{m}+\frac{x_m^{-4}}{m}
\right).
\]
Note that $|a_m| \leq c/m^2$ and $|b_m| \leq c/m$, for some $c > 0$
small enough.

Take
\[
\la_m = \frac{1+Z_m}{3m}.
\]
We have that
\begin{eqnarray*}
\frac{1+Z_{m-1}}{3(m-1)} &=& \frac{2}{9m^2} + a_m +
                           \left(1 + \frac{2}{3m} + b_m\right) \times \nonumber \\
&& \times
                            \left(\frac{1}{3m} + \frac{Z_m}{3m}\right)
     \left(1 - \frac{1}{3m} - \frac{Z_m}{3m} + O\left(\frac{1}{m^2}\right)
                      +O\left(\frac{Z_m^2}{m^2}\right)\right) \nonumber \\
&=& \frac{2}{9m^2} + a_m + X_1 \cdot X_2 \cdot X_3.
\end{eqnarray*}
\begin{eqnarray*}
X_2\cdot X_3 &=& \frac{1}{3m} + \frac{Z_m}{3m} - \frac{1}{9m^2} -
\frac{2Z_m}{9m^2} -  \nonumber \\
&& - \frac{Z_m^2}{9m^2} + O\left(\frac{1}{m^3} +
\frac{Z_m}{m^3} +\frac{Z_m^2}{m^3}+ \frac{Z_m^3}{m^3}\right), \nonumber
\end{eqnarray*}
\begin{eqnarray*}
X_1\cdot X_2 \cdot X_3 &=& \frac{1}{3m} - \frac{1}{9m^2} + \frac{b_m}{3m}
+ \frac{Z_m}{3m} + \frac{b_mZ_m}{3m} - \nonumber \\
&& -
\frac{Z_m^2}{9m^2} +  O\left(\frac{1}{m^3} +
\frac{Z_m}{m^3} +\frac{Z_m^2}{m^3}+ \frac{Z_m^3}{m^3}\right), \nonumber
\end{eqnarray*}
Therefore
\begin{eqnarray*}
Z_{m-1} &=& R_m + Z_m \times \nonumber \\
&&
\left(1-\frac{1}{m} + b_m + O\left(\frac{1}{m^2}\right)-Z_m
\left(\frac{1}{3m}+ O\left(\frac{1}{m^2}\right)\right) +
O\left(\frac{Z_m^2}{m^2}\right)\right), \nonumber
\end{eqnarray*}
where
\[
R_m =  3ma_m + b_m +O\left(\frac{1}{m^2}\right).
\]

If we fix a small $\de>0$, then for $n$ large enough
\[
|Z_{m-1}| \leq |R_m| + |Z_m|\left(1-\frac{\delta}{m}\right).
\]
Without affecting the asymptotic behavior of $Z_m$, we can assume
that the bound above holds for all $m\geq 3$. Using it
recurrently we get
\begin{eqnarray*}
|Z_m| &\leq& \sum_{k=m}^{N-N_3}|R_k|\prod_{i=m}^{k}
                            \left(1-\frac{\de}{i}\right) \nonumber \\
&\leq & \const\sum_{k=m}^{N-N_3}\left(|R_k|\exp\left(-\sum_{i=m}^{k}
                         \frac{\de}{(i)}\right)\right) \nonumber  \\
&\leq & \const\sum_{k=m}^{N-N_3}|R_k|(m/k)^{\de}.
\end{eqnarray*}
Then
\begin{eqnarray*}
\sum_{m=2}^{N-N_3}|\chi_m| &\leq& \sum_{m=2}^{N-N_3}|Z_m|/m
                                                           \nonumber  \\
    &\leq&        \const\sum_{m=2}^{N-N_3}
                        \sum_{k=m}^{N-N_3}|R_k|m^{\de-1}/k^{\de}
                                                          \nonumber   \\
&\leq&         \const\sum_{k=2}^{N-N_3}
                     |R_k| \sum_{m=2}^{k}m^{\de-1}/k^{\de} \nonumber \\
&\leq&           \const\sum_{k=2}^{N}|R_k|.
\end{eqnarray*}

The last sum is uniformly bounded on $N$, which completes the proof.
\end{proof}

\begin{proof}[Proof of Theorem \ref{teo:derivada}]
Let $dx$ be an unstable vector.
At the exiting period, $\lambda_m \asymp 1/m$ and $\tau_m \asymp m^{-2/3}N^{1/6}$,
for $m=2,\ldots,N-N_3$. Therefore  $B_m^{+} = \frac{\lambda_m}{\tau_m} \asymp m^{-1/3}N^{-1/6}$,
$m=2,\ldots,N-N_3$. When $m=1$, $\tau_N \asymp N^{1/6}$.

For $m=1$ (or $n=N$),
\[
 B_N^{+} \asymp B_{N-1}^{+} \asymp N^{-1/6}.
\]

Hence, from Lemma \ref{lema:lambda_1_N1} and Lemma
\ref{lema:lambda_N3_N}, we get
\[
 \frac{\|D_xT^{N_5+1}(dx)\|_p}{\|dx\|_p} \asymp  N^{2/3}\times N^{1/3} \asymp N,
\]
for all $x\in E_N$ satisfying $\ga_1\asymp N^{-1/3}$ and $\ga_N \asymp N^{-1/3}$.
\end{proof}

\section{Proof of Theorem \ref{main:decay}}

\begin{figure}[htpb]
\begin{center}
 \includegraphics[height=7cm,width=9cm]{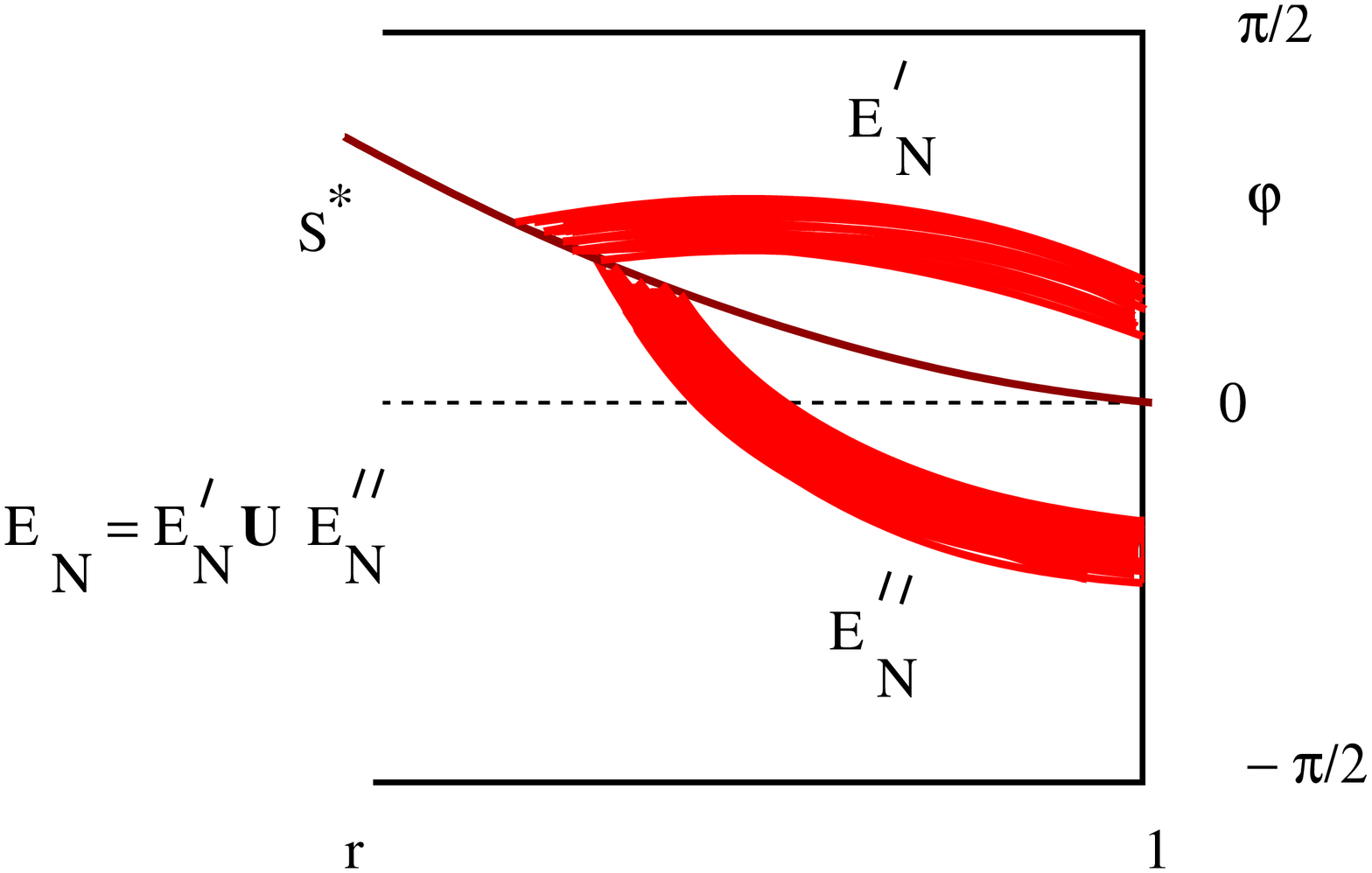}
\caption{\label{fig:EN} The sets $E_N^{\prime}$, $E_N^{\prime\prime}$ and $E_N$.}
\end{center}
\end{figure}

Let $E_N= \{x \in M_4 \;\;|\;\; R(x)= N+1\}$.
 This set is as in Figure \ref{fig:EN}.
It is bounded by curves, denoted by $S^*$, $S_{N-1}^*$ and $S_{N}^*$,
and by the line $r=1$. 
The curve $S^*$ is made up of points from $M_4$ that, leaving $\Le$,
they hit the dispersing part
$\U$ tangentially at the first collision. This is a decreasing curve, 
since it is a singularity line for $T_5$ for positive values of $\varphi$
 until $(1,0) \in M_5$;
 and it is not hard to calculate the slope of this line,
 obtaining that it has an horizontal tangency at  $(1,0) \in M_5$
Besides, the lines $S^*_N$ separating $E_N$ and $E_{N+1}$ 
are constituted of trajectories which the last collision in $\U$ leaving out 
the cusp is tangent. So, they are singularity lines for $T^N_5$, and then ,
decreasing lines and regular as consequence of the results in \cite[Chapter 4]{CM06}.

 The images $F_N = T_4(E_N)=T^{N+1}_5(E_N)$
 are domains bounded by singularity lines for
 $T_5^{-i}$, $i=1,2,\ldots, N $,
 which are curves  with positive slope. Moreover, by the property
of time-reversing of the billiard map, 
$(r,\varphi) \in E_N$ if, and only if,, $(r,-\varphi) \in F_N$. So $F_N$ is 
obtained reflecting $E_N$ along the line $\varphi = 0$.

The domain $E_N$ close to $(1,0) \in M_5$ is made up of two strips:
 the inferior strip $E_N^{\prime\prime}$, which consists of points
 that leave the vertical wall and hit the cusp directly;
 and the superior strip $E_N^{\prime}$,
 which consists of points that hit the cusp after a rebound in the horizontal
 part of the table $\D$. The sets $E_N$, $N>N_0$, make up a nested
 structure that shrink to (1,0) as $N$ goes to infinity (it is enough 
 to note that in order to achieve more rebounds inside the cusp, on $\D$,
 we must begin closer to the point $(0,0)$ and the particle must be thrown
 almost parallel with respect to the axis $x$).

The point from $E_N$ farthest from $(1,0) \in M_5$ over $S^*$ is at a distance
 $\asymp N^{-1/6}$ because $x_1 \asymp N^{1/6}$ on $\D$  (given by equation (\ref{eq:x1})).
 Over $r=1$, using the values of $x_1$ and $\gamma_1$ obtained in $(\ref{eq:x1})$ and 
$(\ref{eqn:gamma1})$, respectively, and a simple geometric construction,
we get that the distance of the strip $E_N^{\prime\prime}$ to the point $(1,0)$ is $\asymp N^{-1/3}$.
Since the lines $S^*$, $S^*_{N-1}$ and $S^*_N$ are decreasing and $S^*$ has 
horizontal tangent on$(1,0)$, the ``length"  of each strip of $E_N$ is $\asymp N^{-1/6}$.

Now consider an unstable curve $W$ inside one of the strips of $E_N$,
transverse to the direction of $S^*_N$, by the relation between cones
and singularity lines given by condition (C6) from \cite[Section 8]{Le02}. 
 Using the symmetry of $T_5^{N+1}(E_N)$, the set 
$T_5^{N+1}(W)$ is a line stretching "from top to bottom" one of the strips of
 $F_N=T_5^{N+1}(E_N)$, therefore, it has ``length" $\;\asymp N^{-1/6}$.
Using the fact that the derivative of $T_5^{N+1}$ has an expansion rate of
 $\asymp N$ for unstable vectors, given by Theorem \ref{teo:derivada},
we get that $|W| \asymp N^{-1/6}/N = N^{-7/6}$.
This is the ``width" of each of the strips of $E_N$.

\begin{figure}[htpb]
\begin{center}
 \includegraphics[height=7cm,width=9cm]{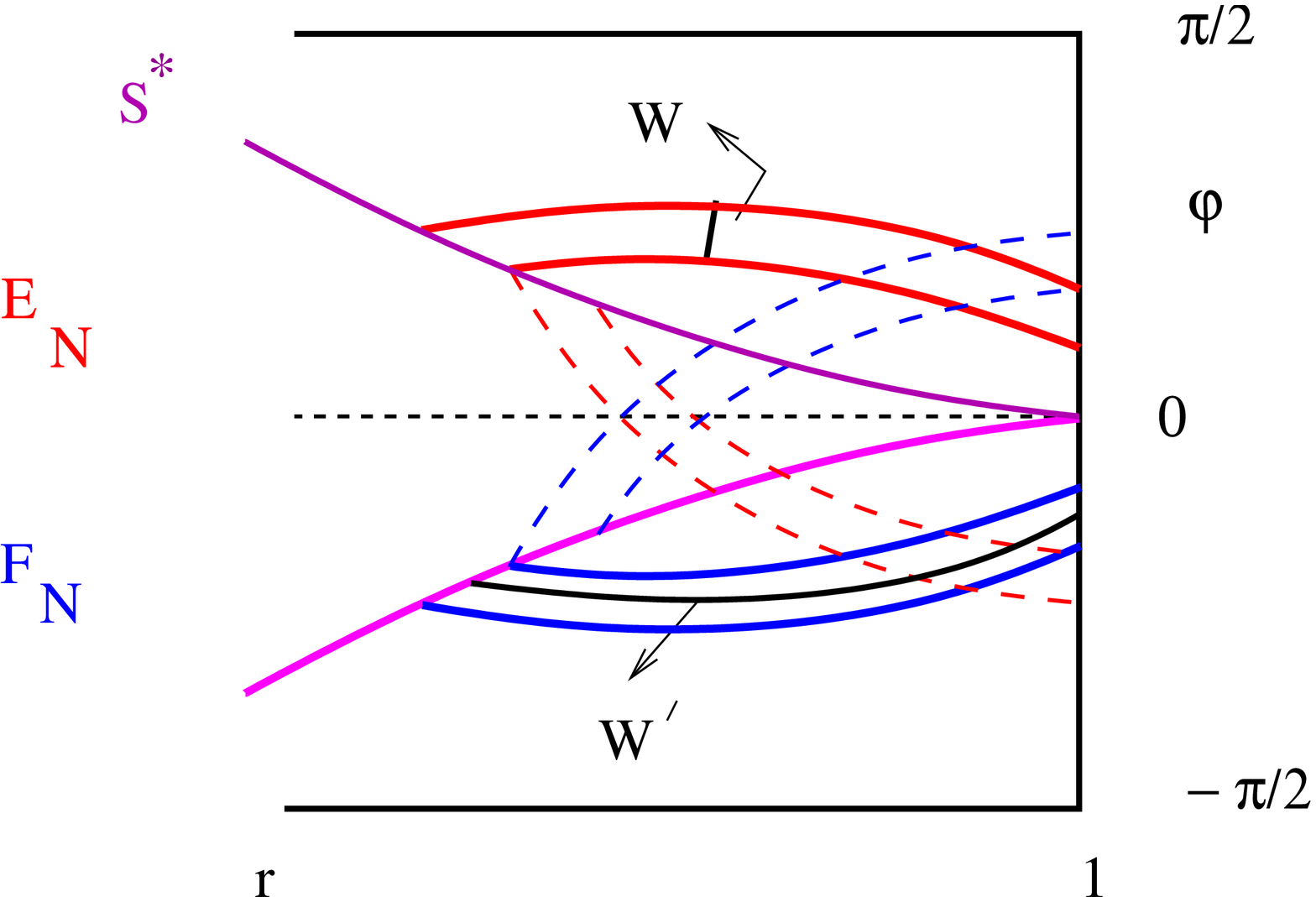}
\end{center}
\end{figure}
Since the sets $E_N$ are away from $\varphi = \pm \pi/2$, the measure $\mu$
is equivalent to the Lebesgue measure on $\R^2$.
Thus,
\[
\mu(E_N) \asymp N^{-1/6} \times N^{-7/6} = N^{-4/3},
\]
so
\[
\mu\left(\bigcup_{n=N}^{\infty}E_n\right) \asymp N^{-1/3},
\]
 hence $A = \bigcup_{n=N}^{\infty}E_n$ has finite measure.

Also the measure of the intersection
 $E_{m} \cap T_5^{m}E_m$ can be computed using the symmetry of the sets
$E_m$ and $F_m$,
\[
 \mu\left(E_{m} \cap T_5^{m}E_m\right) \asymp m^{-7/6} \times m^{-7/6} = m^{-7/3}.
\]
Thus
\[
 \mu\left(A \cap T_5^m A\right) = \mu \left(\bigcup_{n=N}^{\infty}E_n
\cap T_5^m\left(\bigcup_{n=N}^{\infty}E_n\right)\right) \geq
\mu\left(E_{m} \cap T_5^{m}E_m\right) \asymp m^{-7/3},
\]
showing that the speed of decay is at most polynomial.

\section*{Acknowledgments}
This work corresponds to the Ph.D. thesis of the fourth author, under
the guidance of Maria Jose Pacifico and Roberto Markarian. He thanks
IMERL-Montevideo for the hospitality during certain period while
working on this project.

R. Markarian and R. Soares thank UNAM, Instituto de Matematicas
in Cuernavaca, Mexico, and A. Arbieto, R. Markarian and 
R. Soares thank ICTP-Trieste for their kind hospitality and
 financial support during a visit there while working on this project.

M.J. Pacifico thanks Scuola Normale Superiore di Pisa for its kind
hospitality during her staying for a postdoctoral program, when this 
work was finished.


\appendix
\section{Mixing systems and another proof of Corollary
 \ref{main:mixing}} 
\label{sec:mainmixing}

In this section we explain a sufficient condition for a system to be F-mixing.
It is based on the work of Coudene \cite{Co07b}
 adapted to our definition of F-mixing. From now on
$X$ is a metric space, $\mathcal{A}$ the Borel $\sigma$-algebra of $X$, 
$\mu$ an infinite $\sigma$-finite regular measure on $X$ and $\F:
X \rightarrow X$ a $\mu$-measure preserving transformation.


\begin{definition}{\rm (\cite[Definition 1]{Co07b})}
We define the {\it stable distribution} of $\F$ of a point $x \in X$ as
\[
W^s(x) = \{ y \in X : d(\F^n(x),\F^{n}(y)) \to 0\;\; \mbox{as}\;\; n \to \infty\}.
\]

A measurable function $f: X \rightarrow \R$ is called {\it $W^s$-invariant}
when there exists a set $\Omega \subset X$ with full measure such that
for all $x,y \in \Omega$, $y \in W^s(x)$ implies $f(x)=f(y)$.
\end{definition}

If $\F$ is invertible we define the {\it unstable distribution}
$W^u(x)$ of a point $x$ for $\F$ as the stable distribution
for $\F^{-1}$. In a similar way , we define a $W^{u}$-invariant function.

We say that the stable distribution $W^s$ is ergodic if every $W^s$-invariant function is $\mu$-almost
everywhere.

The propositions below are slight modifications of Theorem 2 and Theorem 3
from \cite{Co07b}. Indeed, since the main elements used in the proof are 
the Banach-Alaoglu Theorem and Banach-Saks Theorem, which are true in Hilbert 
spaces, there will be few changes in the proofs. We also assume that the measure is
regular because we use the fact that continuous functions with compact support
are dense in $L^2_{\mu}(X)$ \cite[p.69]{Ru87}.

\begin{proposition}{\rm (Based on \cite[Theorem 2]{Co07b})} \label{prop:wsinv}
Let $X$ be a metric space, $\mu$ a regular infinite $\sigma$-finite measure
on $X$, $\F :X \rightarrow X$ a $\mu$-measure-preserving transformation
and $f \in L^2_{\mu}(X)$.
Then any weak limit of $f \circ \F^n$ is $W^s$-invariant.
\end{proposition}

\begin{proof}
Let $g$ be a weak limit of $f\circ \F^{n_i}$. First assume that
$f$ is continuous with compact support (thus uniformly continuous). The Banach-Saks
theorem guarantees that there exist subsequences $m_l$ and $n_{i_k}$ such that
\[
\Psi_l(x) = \frac{1}{m_l}\sum_{k=1}^{m_l}f \circ \F^{n_{i_k}}
\stackrel{l \to \infty}{\longrightarrow}
 g \;\;\; \mu-q.t.p..
\]

If $y \in W^s(x)$, then
\[
|\Psi_l(x) - \Psi_l(y)| \leq \frac{1}{m_l}\sum_{k=1}^{m_l}|f \circ
\F^{n_{i_k}}(x)-  f \circ \F^{n_{i_k}}(y)| \stackrel{l \to \infty}{\rightarrow} 0.
\]
So $g$ is $W^s$-invariant.

Let $f \in L^2_{\mu}$. For all $\varepsilon > 0$, there exists 
a continuous function $f_0$ with compact support such that
$\|f-f_0\|_2 < \varepsilon$. 
Passing to a subsequence, by Banach-Alaoglu Theorem,
we can assume that $f_0 \circ \F^{n_{i}}$ converges weakly to a 
function $g_0$ 
which is $W^s$-invariant. It follows that
$(f-f_0) \circ \F^{n_{i}} \longrightarrow g-g_0$ weakly, which implies that
\[
\|g-g_0\|_2 \leq \liminf \|(f-f_0) \circ \F^{n_{i}}\|_2 \leq
\|f-f_0\|_2 < \varepsilon.
\]
Thus there exists a sequence of $W^s$-invariant functions that converges 
to $g$ in the $L^2_{\mu}$-norm and, passing to a subsequence,
almost everywhere.
Hence, for a set $\Omega$ with full measure, if $y,x \in \Omega$, 
$y \in W^s(x)$, we get that
\[
g(y) = \lim g_n(y) = \lim g_n(x) = g(x).
\]
 This shows that $g$ is $W^s$-invariant.
\end{proof}

Using the proposition above, the next one is proved
as in \cite{Co07b}.

\begin{proposition} {\rm (Based on \cite[Theorem 3]{Co07b})}
Let $X$ be a metric space, $\mu$ a regular infinite $\sigma$-finite measure
on $X$, $\F :X \rightarrow X$ an invertible $\mu$-measure-preserving transformation
and $f \in L^2_{\mu}(X)$.
Then any weak limit of $f \circ \F^n$ is $W^s$-invariant and $W^u$-invariant.
\end{proposition}

\begin{corollary} \label{cor:wserg}
If $W^s$ is ergodic then $\F$ is F-mixing.
\end{corollary}

\begin{proof}
Let $f \in L^2_{\mu}$. If $f \circ \F^n$ has a weak limit,
by Proposition \ref{prop:wsinv} and hypothesis, this limit
is constant at almost every point; thus it is equal to zero
at almost every point. Therefore $F$ is F-mixing.

Suppose that $f \circ \F^n$ does not converge weakly to zero. 
Thus there exist an $\varepsilon > 0$, a subsequence $n_i$ and a function
$h \in L^2_{\mu}$ such that
\[
\lim_{i \to \infty} \int  (f \circ \F^{n_i})h\; d\mu > 0.
\]

However by Banach-Alaoglu Theorem 
there exists a subsequence $n_{i_k}$ such that $f \circ \F^{n_{i_k}}$ converges
weakly to a function $W^s$-invariant, by
Proposition \ref{prop:wsinv}, which is constant by hypothesis and hence 
must be zero almost everywhere. This contradiction shows that
$\F$ is F-mixing.
\end{proof}

\begin{proof}[Proof of Corollary \ref{main:mixing}]
Take a function $\psi:X \rightarrow \R$ which is
$W^s$-invariant and $W^u$-invariant. By the property of absolute continuity 
of the local stable and unstable manifolds \cite{Le02}, Theorem 7.5,
$W^s$-invariance and $W^u$-invariance imply that this function must be constant
almost everywhere in the ergodic component of $T_5$. However, since $T_5$
has only one ergodic component
 $\psi$ is constant almost everywhere,
that is, $W^s$ and $W^u$ are ergodic. Thus, by Corollary \ref{cor:wserg},
$T_5$ is F-mixing
\end{proof}

 A. Arbieto, M. J. Pacifico, R. Soares:
Instituto de Matem\'atica,
Universidade Federal do Rio de Janeiro, 
C. P. 68.530, CEP 21.945-970, 
Rio de Janeiro, R.J., Brazil.  \\
 E-mails:                     
{\it arbieto@im.ufrj.br, pacifico@im.ufrj.br,} 
{\it castijos@gmail.com} \\

 R. Markarian: 
 Instituto de Matem\'atica y Estad\'{\i}stica (IMERL),
Facultad de Ingenier\'{\i}a, Universidad de la Rep\'ublica,
CC30, CP 11300, Montevideo, Uruguay.\\
 E-mail: {\it roma@fing.edu.uy}

\end{document}